\documentclass[a4paper,10pt]{article}
\usepackage{amssymb,amsbsy,amsmath,amsfonts,amscd} 
\usepackage{mathptm}
\usepackage{graphicx}

\newtheorem{lem}{Lemma}
\newtheorem{theorem}{Theorem}
\newtheorem{prop}{Proposition}
\newtheorem{corol}{Corollary}
\newtheorem{rk}{Remark}
\newtheorem{an}{Ansatz}

\numberwithin{theorem}{section}
\numberwithin{defi}{section}
\numberwithin{corol}{section}
\numberwithin{prop}{section}
\numberwithin{lem}{section}
\numberwithin{equation}{section}
\numberwithin{rk}{section}
\numberwithin{an}{section}

\def\cof{\mathop{\hbox{\rm cof}}}

\def\o{\omega}
\def\O{\Omega_1}
\def\R{\mathbb{R}}
\def\M{\mathbb{M}}
\def\tr{\mathop{\hbox{\rm tr}}}
\def\w{\widetilde}
\def\x{\xi}
\def\z{\zeta}
\def\qed{\hfill$\square$\par\medbreak}
\def\botimes{\mathbin{\bar\otimes}}
\def\nablap{\nabla_{\!\!p}}

\def\det{\mathop{\hbox{\rm det}}}

\def\beq{\begin{equation}}
\def\eeq{\end{equation}}

\def\R{\mathbb {R}}

\def\a{\alpha}
\def\o{\omega}
\def\p{\psi}
\def\P{\Psi}
\def\f{\varphi}

\def\grd{\nabla}
\def\grdc{{\nabla}^2}

\def\nor{\arrowvert}
\def\w{\widetilde}
\def\O{\Omega_1}
\def\dis{\displaystyle}
\def\h{\frac{1}{h}}
\def\hc{\frac{1}{h^2}}
\def\pt{\otimes}
\def\cof{\mathop{\hbox{\rm cof}}}
\def\tr{\mathop{\hbox{\rm tr}}}
\def\z{\zeta}
\def\x{\xi}

\def\cv{\rightarrow}
\def\ptc{\overline{\otimes}}

\def\r{\frac{1}{r}}

\begin{document}
\title{Asymptotic analysis for a second order curved thin film}


\author{Hamdi Zorgati\thanks{Imam Mohammad Ibn Saud Islamic University (IMSIU), College of Sciences, Department of Mathematics and Statistics, PO-Box 90950, Riyadh 11623 Saudi Arabia and University of Tunis El Manar, Faculty of Sciences of Tunis, Department of Mathematics, Tunisia.
		E-mail addresse:  hzorgati@imamu.edu.sa }}

\date{}
\maketitle

\begin{abstract}
We consider a second order thin curved film whose behavior is
governed by an energy made up of a first order nonlinear part depending on the gradient of the deformation augmented by a
quadratic second order part depending on the tensor of
second derivatives of the deformation. We carry out a 3D-2D analysis through an
asymptotic expansion in powers of the thickness of the film as it tends to zero.


\medskip

\end{abstract}

\section{Introduction}

We are interested in this study in some second order functionals that model the mechanics of thin curved films. We study the behavior of
such functionals as the thickness of the thin film goes to zero by performing a 3D-2D analysis. The energies considered contain a first order part depending on the gradient of the deformation, as in the case of hyperelastic materials, through a nonlinear density. The second part of the energy contains a quadratic  term depending on the tensor of the second derivatives of the deformation. When the quadratic form is the square of the Euclidean norm, we recover the example of the interfacial energy term of Van der Waals type for martensitic materials (see \cite{bhajam,ledzor,ledzor2,shu,ledmeu,ledmeu2}). 

The 3D-2D dimension reduction for thin films is an efficient way to obtain a two-dimensional model that approximates a three-dimensional realistic one when the thickness is small enough. Two-dimensional approximations are easier to study and simpler in terms of numerical simulations with less cost in data and time solving (see \cite{Ber,chb}).

Basically, there are two main ways in proceeding with the dimension reduction. The first method, that we consider in this work, is the formal asymptotic expansion method. It is based on an Ansatz assuming that the deformation of the three-dimensional thin film admits an asymptotic expansion in powers of the thickness in order to characterize the limit behavior.
The method of formal asymptotic expansion for plates in a functional framework was introduced by P.G. Ciarlet
and P. Destuynder in \cite{ciades} and P. G. Ciarlet \cite{cia0}, after the pioneering works of K. O. Friedricks and R. F. Dressler \cite{FD} and A. L. Goldenveizer \cite{Gol}. Many authors used this method in performing dimension reduction for different models such as \cite{SP1,SP2,SP3,cialod,cialodmia,mia,CR1,CR2,CDP}. Another approach of this method using variational formulation only
in terms of displacements was derived by A. Raoult in \cite{rao}.
A general description of this method and of works where
it was used can be found in P. G. Ciarlet \cite{cia31} (see also \cite{mia,cia2,cia3,cia4}). Later on, this method was
refined and improved in term of assumptions by D. D. Fox, A. Raoult and
J. C. Simo in \cite{foxraosim}. Their method is based on the resolution of a sequence of Euler-Lagrange
equations.  

The second approach for dimension reduction is variational. It needs no ansatz on the deformations and uses arguments of $\Gamma $-convergence, a kind of convergence that was introduced by E. De Giorgi in \cite{DeG, DF}. The main properties of the $\Gamma$-convergence are first that, up to a subsequence, the $\Gamma$-limit always exist and second that, if a sequence of almost minimizers stays in a compact then the limits of the converging subsequences minimize the $\Gamma$-limit. Thus, the limit model obtained, provides a description of the behavior of the minimizers of the sequence of energies depending on the thickness, as is goes to zero. Also, the limit obtained being lower semi-continuous with respect to the considered topology, the existence of minimizers  is thus guaranteed (see \cite{DAL} ). This method was used in many works in the context of dimension reduction such as \cite{ABP,BFF,FF,IZ,ledrao,ledrao2,Zor1,Zor2,Zor3,bZZ}.  

In this work, we use the formal asymptotic expansion method introduced by O.~Pantz in \cite{pan} which leads to the resolution of a succession of
minimization problems (see \cite{Tra,Meu}). Indeed, it appeared that the formal expansion method yields an incorrect result compared to the $\Gamma$-convergence method in some cases. Pantz's method makes it possible to obtain formal expansions that agree with $\Gamma$-convergence limits as much as is possible   

The equilibrium state of the thin curved film 
occupying the domain $\widetilde{\Omega}_{h}\subset \mathbb{R}^3$ of thickness $h$ with midsurface a bounded two-dimensional submanifold $\widetilde{S}$, and undergoing a
deformation $\w \f$ is thus described by the minimizers of the energy
\begin{equation}
\w{J}^h (\w{\f})  = \int_{\widetilde{\Omega}_{h}}
\big[\, Q(\grdc\widetilde{\f} )\,+\,W ({\grd
  \widetilde{\f}})\big] \, dx 
\end{equation}
where $\grdc\widetilde{\f}$
is the $3\times 3 \times 3$ tensor of second derivatives, $Q$ is a positive  quadratic form and $
W({\grd \widetilde{\f}})$ represents the elastic energy by
unit volume of the material depending on the gradient of
the deformation $\widetilde{\f}$. We are interested in the
behavior of this energy and its minimizers when the thickness $h$ of
the film goes to zero. 

We begin our study by introducing the notation and some geometrical
preliminaries, then we rescale the energy in order to work on a domain
independent of the thickness $h$. Next, we carry out the computations of the terms appearing in the asymptotic expansion
of the energy, which will bring us to the resolution of the minimization problems.
In order to be able to perform such an analysis, we will consider a first order nonlinear term given by the Saint Venant-Kirchhoff material stored energy function.
 
 The term of order zero obtained by solving the first four minimization problems  will provide the two-dimensional scaled limit energy that we illustrate with an example in the case of
plates and also in some simple examples of curved films. Notice that we recover, as a particular case supposing the quadratic form to be the square of the norm, the two-dimensional second order energy part, obtained by convergence arguments in \cite{bhajam} for the planar films and \cite{ledzor,ledzor2} for the curved films. Moreover, the first order energy part generalizes to the case of curved films the results obtained in \cite{pan}.
  

\section{Notations and geometrical preliminaries}
Let $(e_1, e_2, e_3)$ be the orthonormal canonical basis of the
Euclidean space $\mathbb{R}^3$. The norm of a vector
$v$ of $\mathbb{R}^3$ is denoted by $\nor v \nor$, the scalar
product of two vectors $u$ and $v$ by $u\cdot v$, their vector
product by $ u\wedge v$ and their tensor product by $u
\otimes v$. Let $\M_{33}$ be the space of $3\times 3$ matrices with 
real coefficients equipped with the usual norm $\nor F
\nor =\sqrt{ \tr (F^T F)}$. We denote by $A = (a_1\nor a_2\nor
a_3)$ the matrix of $\M_{33}$ whose first column is $a_1$ the
second is $a_2$ and the third is $a_3$. We denote by $\grd \f$ and
$\grdc{\f } $ the $3\times 3$ and $3\times 3\times 3$ tensor of first and second derivatives
of a deformation $\f:\Omega\cv \R^3$, where $\Omega$ is an open subset of $\mathbb{R}^3$. 

Let $ \mathbb M _{333}$ be the space of $3\times 3\times 3$ 
tensors. When $ P = (p_{ijk }) $ is a tensor of
order 3, the tensor with components $r_{ijk } = p_{ikj}$ is the transpose
of $P$ with respect to the second and third index. Let $P$ be a tensor
of order $p\geq 1$ and $R$ a tensor of order $r\geq 1$. We
define the contracted tensor product of the tensors $P$ and $R$ and
we denote by $P\ptc R$, the tensor of order $p+r-2$ whose components are
 obtained by summation of the products of the components of $P$ and
$R$, on the last index of the components of $P$ and the first index of
those of $R$. For example, using the summation convention, if $p=3$ and $r=2$ we have
$$   
P \ptc R = P_{ijk} R_{km} (e_i\pt e_j \pt e_m).
$$
Let $Q: \mathbb M _{333} \cv \mathbb R^+$ be a positive quadratic form whose corresponding bilinear symmetric form is $B : \mathbb M _{333} \times \mathbb M _{333} \cv \mathbb R$. In order to simplify our computations, we will make the following assumptions on the quadratic form $Q$ : We suppose that the basis $(e_i\pt e_j \pt e_k)_{1\leq i,j,k\leq 3}$ of $\mathbb M _{333}$ is orthogonal with respect to the quadratic form $Q$, i.e. the matrix of the quadratic form in this basis is diagonal so that for every $A\in \mathbb M _{333}$, we have 
$$
Q(A)= \sum_{1\leq i,j,k\leq 3} A_{ijk}^2 Q_{ijk},
$$ 
where $A_{ijk}$ denote the components of the tensor $A$ in the orthogonal basis and $Q_{ijk}= Q(e_i\pt e_j \pt e_k)$. Notice that some components $Q_{ijk}$ might vanish since we do not suppose $Q$ to be definite. Nevertheless, in order to solve the minimization problems, we suppose that at least one of the $Q_{ijk}$ is strictly positive.

Throughout this work, we will use the convention of summation. The
Greek indices take their values in the set $ \{ 1, 2\}$ and the Latin
indices take their values in $\{1, 2, 3\}$.

We consider an open domain $\widetilde{\Omega}_{h}$ occupied by
a thin curved film of thickness $h$ whose midsurface
is $\widetilde{S}$, a two-dimensional bounded submanifold of class $
C^2$ of $ \mathbb{R}^3$ that admits an atlas
consisting in one chart. Let $\p$ be this chart.  It
is thus a $ C^2$-diffeomorphism from an open bounded domain $\omega\in
\mathbb{R}^2$ into $ \widetilde{S} $. We
suppose that $ \omega$ has a Lipschitz boundary $
\partial{\omega}$ and that $\p$ is extendable into
a function of $ C^2(\overline{\omega}, \mathbb{R}^3)$. Let $
a_{\alpha}(x) = { \p_{, \alpha}}(x)$ be the vectors of the covariant basis
of the tangent plane $T_{\p (X) } \widetilde{S}$
associated to the chart $ \p$, where $ { \p_{, \alpha}}$ denote
the derivatives of $ \p$ with respect to the variable
$\alpha$. We suppose that there exist $\delta > 0$ such that
$$
\nor a_1 (x)\wedge a_2 (x) \nor \geq \delta\;\text{ on}\; \overline{\omega},
$$
and we define the vector $a_3 (x) = \frac{a_1 (x)\wedge
    a_2 (x)}{\nor a_1 (x)\wedge a_2 (x) \nor } $ belonging to
$ C^1(\overline{\omega}, \mathbb{R}^3)$. This is the normal unit
vector to the tangent plane. The vectors $a_1 (x)$, $a_2 (x)$ and
$a_3 (x)$ constitute the covariant basis. We define the vectors of
the contravariant basis by the relations
$$
 a^{\alpha} (x) \in T_{\p (x)} \widetilde{S},\;a^{\alpha} (x) \cdot
 a_{\beta} (x) = \delta_{\alpha \beta}\;\text{and}\;a^3 (x) = a_3
 (x). 
$$
We set $A (x) = (a_1 (x) \nor a_2 (x) \nor a_3 (x) ) $. We notice that
$ A (x)$ is an invertible matrix on $ \overline{\omega}$ and its
inverse is given by
$$
 A^{-1} (x) = {(a^1 (x) \nor a^2 (x) \nor a^3 (x) )}^T .
$$
We also notice that
$$
 \det A (x) = \nor \cof A (x) \cdot e_3\nor = \nor a_1 (x)\wedge a_2
 (x) \nor > \delta > 0\;\text{on}\; \overline{\omega}. 
$$
The domain $ \widetilde{\Omega}_{h}$ is defined as
$$ 
  \widetilde{\Omega}_{h} = \big\{ x \in \mathbb{R}^3, \exists\, \w{x}
  \in \widetilde{S} ,  x = \w{x} + \eta\,  a_3
  \big(\p^{-1}(\w{x})\big)\;  \hbox {avec} \;  \frac{-h}{2} < \eta <
  \frac{h}{2} \big\} 
$$
and for $\widetilde{x}\in \widetilde{\Omega}_{h}$ we have
$$
\w{x} = \w{\pi} (\w{x}) + \Big[\big(\w{x}-\w{\pi}(\w{x})\big)\cdot
a_3\big(\p^{-1}(\w{\pi}(\w{x}))\big)\Big]\cdot
a_3(\p^{-1}(\w{\pi}(\w{x}))), 
$$ 
where $\w{\pi}$ is the orthogonal projection which sends $
\w{\Omega}_h$ onto $\w{S}$. The map $ \w{\pi}$ is well
defined and is of class $C^1$ for $h$ small enough compared to the
radius of curvature of $\w S$. We can thus define the curvilinear
coordinates of a point of the domain $\w{\Omega}_h$ associated to
the chart $\p$ by
$$
 (x_1, x_2 ) = \p^{-1}(\w{\pi}(\w{x}))\;\text{ and }\;x_3 =
 (\w{x}-\w{\pi}(\w{x}))\cdot a_3(\p^{-1}(\w{\pi}(\w{x}))). 
$$
We now define the domain $\Omega_h$ by
$$
\Omega_h = \big\{ x\in \R^3 , x = (x_1, x_2, x_3 ), (x_1, x_2 )\in
\omega, \frac{-h}{2}< x_3 <\frac{h}{2}\big\}. 
$$
We choose a small enough $h^*$ and define a $C^2$-diffeomorphism $\P$:
$\Omega_{h^*}\rightarrow\w{\Omega}_{h^*}$ by 
\beq
 \P (x_1, x_2, x_3 ) = \p (x_1, x_2 ) + x_3 a_3(x_1, x_2 ). 
\eeq
This diffeomorphism is the inverse of the rescaling that we will
use in order to work on a flat domain, and its gradient verifies
$$
 \grd\P (x_1, x_2, x_3 ) = A (x_1,x_2) + x_3 (a_{3,1}(x_1, x_2 )\nor
 a_{3,2}(x_1, x_2 )\nor 0). 
$$
The matrix  $\grd\P (x_1, x_2, x_3 )$ is everywhere invertible and its 
determinant is the Jacobian of the rescaling. In the following, ${h}$
will denote some sequence of real numbers less than 
$h^*$ and decreasing into zero.
\section{The three-dimensional and rescaled problems}
The thin film occupies, in its reference configuration, the domain $
\w{\Omega}_h$. We associate to each deformation $\w{\f}:\w{\Omega}_h\rightarrow \mathbb{R}^3$ the energy
$$
\w{J}^h (\w{\f}) = \w{K}^h (\w{\f}) + \w{I}^h (\w{\f}), 
$$
where
$$
\w{K}^h (\w{\f}) = \int_{\w{\Omega}_h} Q(\grdc\w{\f}) \,dx
$$
represents the second order energy and 
$$
{\w{I}^h (\w{\f}) =  \int_{\widetilde{\Omega}_{h}}W ({\grd \widetilde{\f}})\, dx} 
$$
is the elastic energy, with $  W : \M_{33} \rightarrow [ 0,
+\infty [ $ regular and verifying the following growth and coercivity
assumptions: There exists $c_1, c_2 > 0$ such that
$$
\forall A\in \M_{33},\; c_1 (\nor A\nor ^2 - 1)\le W (A) \le c_2 (\nor
A\nor ^q + 1)\quad \text{with}\;2\le q<6.
$$

For simplicity, we will only consider homogeneous boundary conditions
of place imposed on the lateral boundary of the film. We thus
introduce the space of admissible deformations 
$$
\widetilde V_h=\bigl\{ \w\phi \in H^2(\widetilde\Omega_h ;\R^3);
\w\phi(\w x) = \w A\w x \text{ on }  \w\Gamma_h \bigr \},
$$
where $\w{A} = (\w{a_1} |\w{a_2} | \w{a_3})\in  \M_{33}$ is a given
constant matrix and
$\w{\Gamma}_h=\Psi(\partial\omega\times[-h/2,h/2])$ is the lateral
surface of $\widetilde{\Omega}_{h}$. Note that due to the growth
condition satisfied by $W$ and the Sobolev embedding Theorem, the
energy functional $\widetilde J_h$ is well defined and takes its
values in $\R$ for $\widetilde\phi\in \widetilde V_h$.  

The minimization problem consists in finding $\w{\varphi}_h\in
\widetilde V_h$ such that 
\begin{equation}\label{ast}
\w{J}_h(\w\varphi_h)= \min_{\w\phi\in\widetilde V_h} \w{J}_h(\w\phi).
\end{equation}
It is fairly obvious that such a minimizer exists if we suppose the quadratic form to be definite, using
weak convergence in $H^2$ that entails strong convergence in $W^{1,q}$ for
$q<6$. This is not the case for general positive quadratic form.

To study the behavior of this energy and its eventual minimizers, we begin by
flattening and rescaling the minimizing problem in order to work on a
fixed cylindrical domain. For clarity, we proceed in two steps,
flatten first and then rescale.  

 If $\w{\phi}_h$ is a deformation of the film in its reference
 configuration, we define $\phi_h\colon\Omega _h\rightarrow \R^3$ by
 setting $\phi_h(x) = \w\phi_h(\Psi (x))$ for all
 $x\in\Omega_h$. Since $\Psi$ is a $C^2$-diffeomorphism, $\phi_h$ is
 in $H^2$ whenever $\w\phi_h$ is in $H^2$. For every such deformation
 $ \phi_h$, we thus set $J_h (\phi_h) = \w{J}_h (\phi_h\circ
 \Psi^{-1})$ and we obtain that  
$$
\begin{aligned}
J_h (\phi_h) &=\int_{\Omega _h}  \Big\{ Q\Big( \Big (\nabla^2
\phi_h\botimes \big(\nabla {\Psi^{-1}} \circ \Psi \big)\Big)^T
\botimes \big(\nabla {\Psi^{-1}} \circ \Psi \big) + \nabla
\phi_h\botimes \big(\nabla^2 \Psi^{-1} \circ \Psi \big) \Big) \\ 
&\qquad\qquad+W (\nabla \phi_h   (\nabla \Psi^{-1}\circ \Psi))\Big\}
\,\det \nabla \Psi \, dx.\\ 
\end{aligned}
$$
This completes the flattening step.

Let us now turn to the rescaling step. Since the total
energy in the limit membrane deformation regime is of order $h$, we
are actually interested in the asymptotic behavior of the energy per
unit thickness $ \frac{1}{h} {J}_h $.  

We define $ z_{h}:\Omega _1\rightarrow \Omega_h$ by
$$ z_{h}(x_1, x_2, x_3 )=(x_1, x_2, h{x_3}).$$
  To each deformation $\phi_h$ on $ \Omega_h$, we associate the
  deformation $\phi(h)$ on $ \Omega_1$ defined by $\phi(h)(x) =
  \phi_h (z_{h}(x))$ and its rescaled deformation gradient 
$$
\nabla_h\phi(h)=\Big(\phi(h)_{,1}\Big| \phi(h)_{,2}\Big| \frac{1}{h}
\phi(h)_{,3}\Big) 
=\nablap \phi(h) + \frac{1}{h} \phi(h)_{,3} \otimes e_3,
$$
and rescaled tensor of second derivatives 
$$
\nabla^2_h\phi(h)= \nablap^2 \phi(h)  + \frac{1}{h} \big(\nablap
\phi(h)_{,3}\otimes e_3+(\nablap \phi(h)_{,3}\otimes
e_3)^T\big)+\frac{1}{h^2} \phi(h)_{,33} \otimes e_3 \otimes e_3, 
$$
where the notations $ \nablap \phi = \phi_{,{\alpha}} \otimes
e_{\alpha}$ and $ \nablap^2 \phi =  \phi_{,\alpha \beta} \otimes
e_{\alpha} \otimes e_{\beta} $ stand for in-plane first and 
second gradients. Note that $\nablap(\phi_{,3})=(\nablap\phi)_{,3}$,
hence the unambiguous notation $\nablap\phi_{,3}$. Note also that $
\nablap \phi$ is $\M_{33}$-valued and $ \nablap^2 \phi$ is
$\M_{333}$-valued. 
   
Finally, we set $\displaystyle J(h)(\phi(h)) = \frac1h J_h(\phi_h)$. Thus, we obtain
a rescaled energy of the form 
$$
J(h)(\phi(h)) = K(h)(\phi(h)) + I(h)(\phi(h)),
$$
with
$$
 I(h)(\phi)=\int_{\O}W(\nabla_h\phi \,A_h)d_h \,dx
$$
and
$$
K(h)(\phi)=\int_{\O}Q\big ( \big(\nabla^2_h \phi\botimes A_h
\big)^T \botimes A_h+ {\nabla}_h \phi 
 \botimes B_h \big) d_h\, dx,
$$
where 
$$ 
d_h (x) = \det \nabla \Psi (z_h(x)),\; A_h (x) = \nabla \Psi ^{-1}(
\Psi(z_h(x)))\text{ and } B_h (x) = \nabla^2\Psi ^{-1}(
\Psi(z_h(x))). 
$$
Problem (\ref{ast}) becomes: Find $\varphi(h)\in V_h$ such that 
\begin{equation}\label{ast1}
J(h)(\varphi(h))= \min _{\phi\in V_h} J(h)(\phi) ,
\end{equation}
with $V_{h}=\{\phi \in H^2 (\O ;  \R^3) ; \phi(x) = \w{A}\Psi (z_h(x)
) \text{ on } \partial \o \times  {]-} \frac{1}{2}, \frac{1}{2}[\}$. 

Since problems (\ref{ast}) and (\ref{ast1}) are equivalent, and that
the existence of solutions to problem (\ref{ast}) is obvious, there
exists a solution $\varphi(h)\in V_h$ for problem (\ref{ast1}) for every
$0<h<h^*$.

We can now compute the formal asymptotic expansion on $\O$. 


\section{Asymptotic expansion}
To be able to carry out the asymptotic analysis, we need an
explicit expression of the internal elastic energy  density $W$.
We thus consider an elastic energy part that corresponds to a Saint
Venant-Kirchhoff material, i.e.  
$$
 W (A)= \frac{\lambda}{8}(\tr (A^T A - I))^2 + \frac{\mu}{4}\tr ((A^T
 A - I)^2) \; \text{for every}\; A \in \M_{33},
$$
where $\lambda > 0$ and $\mu>0$ are the Lam\'e constants of the material.
This choice can seem not adapted to the context of our study, since
the material of Saint
Venant-Kirchhoff does not present phase changes,
whereas we add a second order energy term to penalize phase changes. 
In fact, we choose the Saint
Venant-Kirchhoff material elastic energy just because it allows for relatively
simple explicit calculations.

We denote by $P(h)$ the problem of finding $\f(h)$ verifying :
$$
J(h)(\f(h)) = \inf \{ J(h)(\f), \f\in V_{h} \}
$$
with 
\beq
V_{h}=\big\{\f \in H^2 (\O ;  \mathbb{R}^3) ; \f(x) = \w{A}\circ \P
(x_1, x_2, hx_3 ) \;\text{ on}\;\partial \omega \times  (-
\frac{1}{2}, \frac{1}{2})\big\}, 
\eeq
where $\w{A} = (\w{a_1}\nor\w{a_2}\nor \w{a_3})$ is a constant matrix
of $\M_{33}$.

\begin{an}\label{marial}
We suppose that the solution $\f(h)$ of the problem $P(h)$ admits a formal
expansion in powers of $h$, i.e. 
$$
\f(h) = \sum_{n\geq 0} h^n \f^n
$$
\end{an}
\begin{rk}\label{quickli}
The first and second derivatives of $\P^{-1}$ taken at the point
$\P(x_1, x_2, hx_3)$ admit each one an 
asymptotic expansion in powers of $h$. This can be
obtained using Taylor-Lagrange formula and supposing that
$\p^{-1}$ is sufficiently regular.
\end{rk}
Using Ansatz \ref{marial} and Remark \ref{quickli}, we
will show that the energy $J(h)$ admits an asymptotic expansion in
powers of $h$. We begin by proving the following proposition.
\begin{prop}\label{prop1}
The energy $J(h)$ admits a first formal asymptotic expansion in powers of $h$ 
$$
J(h)(\f)= \sum_{n=-4}^{\infty} h^n J^n_h (\f).
$$
with 
$$
J^{-4}_h(\f) = \int_{\O} \Big( Q( d_{ij}^{-2}\pt e_i\pt e_j) + C^{ijkl}
E_{ij}^{-2} E_{kl}^{-2}\Big) c_0 \,dx, 
$$
\begin{multline*}
J^{-3}_h(\f) =  \sum _{p=-4}^{-3}\int_{\O} \Big( B(d_{ij}^{p+2}\pt e_i\pt e_j,d_{ij}^{-p-5}\pt e_i\pt e_j ) + C^{ijkl} E_{ij}^{p+2}E_{kl}^{-p-5} \Big)c_0 \,\\
+\Big( Q(d_{ij}^{-2} \pt e_i\pt e_j)  +C^{ijkl} E_{ij}^{-2}
E_{kl}^{-2}\Big) c_1 \,dx  
\end{multline*}
and for every $n\geq -2$
\begin{multline*}
J^n_h (\f) = \sum _{p=-4}^n \int_{\O} \Big(B(d_{ij}^{p+2}\pt e_i\pt e_j, d_{ij}^{n-p-2}\pt e_i\pt e_j) +C^{ijkl} E_{ij}^{p+2}E_{kl}^{n-p-2}\Big) c_0 \,dx \\+ \sum _{p=-4}^{n-1}\int_{\O} \Big(B(d_{ij}^{p+2}\pt e_i\pt e_j, d_{ij}^{n-p-3}\pt e_i\pt e_j) 
+C^{ijkl} E_{ij}^{p+2} E_{kl}^{n-p-3}\Big) c_1 \,dx\\
 + \sum_{p=-4}^{n-2}\int_{\O}\Big( B(d_{ij}^{p+2}\pt e_i\pt e_j, d_{ij}^{n-p-4}\pt e_i\pt e_j)
+C^{ijkl} E_{ij}^{p+2} E_{kl}^{n-p-4}\Big) c_2 \,dx, 
\end{multline*}
where $d_{ij}^p$, $ E_{ij}^{p}$ will be defined in
(\ref{p1})(\ref{p2}) and (\ref{d1}), $c_0$, $c_1$, $c_2$ in
(\ref{jac}), $C^{ijkl}$ are the components of the elasticity tensor and $B$ is the symmetric bilinear form associated with the quadratic form $Q$.
\end{prop}
\begin{rk}
We notice that the terms $J^n_h$ obtained
depend on $h$. It is not an asymptotic expansion in
the usual sense of the term. Nevertheless, Remark \ref{quickli}
will enable us to deduce (see Proposition \ref{Propzah}) that $ J(h)$
admits a second asymptotic expansion in powers of $h$ containing
terms independent of $h$ of the form
$$
J(h)(\f)=\sum_{n\geq-4} J^n(\f)h^n.
$$
Nonetheless, the useful expansion is actually the one given by Proposition \ref{prop1}.
\end{rk}
Next, we use Pantz's method for the asymptotic expansion consisting in solving a
succession of minimization problems and which is based on the
following Proposition \cite{pan}. We provide the proof for the convenience of the reader. In the following, we will
indicate by $ \f$ the sequence $ (\f^i)_{i\in \mathbb N}$
and we will write
\beq
J(h)(\f)= J(h)(\f(h)).
\eeq
\begin{prop}\label{Pantz}
The solution $\f = \Big(\f^0, \f^1, \f^2,...\Big)$ of the sequence of
problems $P(h)$ verifies
$$ 
 \f\in \bigcap _{n=-4}^{\infty} Q_n,
$$
with 
$$
  Q_{n+1} =\big\{ \f \in Q_n, J^n (\f) = \inf _{\f'\in Q_n} J^n(\f') \big\}
$$
and
$$
Q_{-4} = \big\{ \f \in C(\overline{\O} ; \mathbb{R}^3)\, , \sum _n
\f^n h^n \in V_h \big\} . 
$$
\end{prop}
\noindent{\bf Proof}
Let $\f(h)=(\f^i)_{i\in \mathbb N}$ verifying 
$$
J(h)(\f(h))=\inf_{\f\in V_h} J(h)(\f),
$$
with 
$$
J(h)(\f)=\sum_{n\geq-4} J^n(\f)h^n.
$$
We will prove by induction, that $\f(h)\in \bigcap _{n=-4}^{\infty} Q_n$. We have already seen that $\f(h)\in Q_{-4}=V_h$. Suppose that $\f(h)\in \bigcap _{k=-4}^{n} Q_k$ and let us prove that $\f(h)\in Q_{n+1}$. By definition of $Q_n$, we have that for every $\tilde{\f}\in Q_n$
\begin{multline*}
J(h)(\tilde{\f})= \sum_{k=-4}^{n-1} h^k \inf_{\psi\in Q_k}J^k(\psi) + h^n J^n(\tilde{\f}) + \sum_{k=n+1}^{\infty} h^k J^k(\tilde{\f})\\
= C_n(h)+ h^n [ J^n(\tilde{\f}) + \sum_{k=n+1}^{\infty} h^{k-n} J^k(\tilde{\f})].	 
\end{multline*}
Since $\f(h)$ is a minimizer of $J(h)$ over $V_h \supset Q_n$, we have 
$$
J(h)(\f(h)) = \inf_{\psi\in Q_n} J(h)(\psi).
$$
Thus, for every $\psi \in Q_n$, we have 
$$
J^n(\f(h)) + \sum_{k=n+1}^{\infty} h^{k-n} J^k(\f(h))\leq J^n(\psi) + \sum_{k=n+1}^{\infty} h^{k-n} J^k(\psi).
$$ 
Next, for every $\alpha > 0$, there exist $\psi_\alpha\in Q_n$ such that
$$
J^n(\psi_\alpha) \leq \inf_{\psi\in Q_n} J^n(\psi) +\alpha.
$$
Thus, we obtain that
\begin{multline*}
	J^n(\f(h)) + h \sum_{k=0}^{\infty} h^{k} J^{n+k+1}(\f(h))\leq J^n(\psi_\alpha) + h \sum_{k=0}^{\infty} h^{k} J^{n+k+1}(\psi_\alpha)\\
	\leq \inf_{\psi\in Q_n} J^n(\psi) +\alpha + h \sum_{k=0}^{\infty} h^{k} J^{n+k+1}(\psi_\alpha).
\end{multline*}
Letting $h\rightarrow 0 $, we obtain 
$$
J^n(\f(h))\leq \inf_{\psi\in Q_n} J^n(\psi) +\alpha
$$
and thus the result by letting $\alpha \rightarrow 0$.
\qed
We denote by $P_n$ the problem consisting in finding the minimizers of
$J^n$ on $Q_n$. 

Proposition \ref{Pantz} enables us to obtain the solution $\f(h)$ of
$ P(h)$ solving the sequence of problems $P_n$. 
\subsection{Proof of Proposition \ref{prop1} and computation of
  the asymptotic expansion terms} 
Since the first order energy term is of a Saint
Venant-Kirchhoff material type, we have
$$ 
\w{I}^h (\w{\f})= \int _{\w{\Omega }^h} W (\grd \w{\f}),
$$
with
$$
 W (\grd \w{\f}) = C^{ijkl} \w{E}_{ij}^h (\w{\f})\w{E}_{kl}^h (\w{\f}) ,
$$
where
$$
 C^{ijkl} = \lambda \delta ^{ij} \delta ^{kl} + \mu (\delta ^{ik}
 \delta ^{jl} + \delta ^{il} \delta ^{jk}), 
$$
is the elasticity tensor and
$$
 \w{E}^h (\w{\f}) = \frac{1}{2} [\grd \w{\f}^T\grd \w{\f} - I]
$$
is the strain tensor.
\begin{lem}\label{l1}
The Jacobian of the rescaling process $\det \grd \P (x_1, x_2, hx_3)$ is
polynomial with respect to $h$ with
\begin{equation}\label{detz}
\det \grd \P (x_1, x_2, hx_3) = c_0 + h c_1+ h^2 c_2.
\end{equation}
\end{lem}
\noindent{\bf Proof}
Expanding the determinant with respect to the
first and second column of the gradient, we obtain
\begin{equation}\label{jac}
 c_0 = \det (a_1\nor a_2\nor a_3 ),\;  c_1 = x_3[\det (a_1\nor
 a_{3,2}\nor a_3 ) + \det (a_{3,1}\nor a_2\nor a_3 )]\; \text{and}\;
 c_2 = x_3 ^2 \det (a_{3,1}\nor a_{3,2}\nor a_3 ). 
\end{equation}
\qed
\begin{lem}\label{toy2}
The components $ \w{E}_{ij}^h (\w{\f})$ of the strain tensor $
\w{E}^h (\w{\f})$ admit an asymptotic expansion in powers of $h$,
$$ 
 \w{E}_{ij}^h (\w{\f}) = \sum _{n\geq -2} \w{E}_{ij}^n h^n,
$$
where
\begin{equation}\label{p1} 
\w{E}_{ij}^{-2} = \frac{1}{2}\w  C_{ij}^0,\,\w{E}_{ij}^{-1} =
\frac{1}{2}(\w B_{ij}^0 +\w  C_{ij}^1),\,  
 \w{E}_{ij}^{0} = \frac{1}{2}(\w A_{ij}^0 + \w B_{ij}^1 + \w C_{ij}^2
 - \delta _{ij})
\end{equation} 
and 
\begin{equation}\label{p2}
\w{E}_{ij}^{n} = \frac{1}{2}(\w A_{ij}^n +\w  B_{ij}^{n+1} + \w
C_{ij}^{n+2} )\, pour \, n > 0, 
\end{equation}
with 
\begin{equation}\label{p3}
\w A_{ij}^n =\P _{\a,i}^{-1} \P _{\beta,j}^{-1} \sum _{p=0}^n
\f_{,\a}^p \cdot \f_{,\beta}^{n-p}, 
\end{equation}
\begin{equation}\label{p4}
 \w B_{ij}^n = \Big(  \P _{3,i}^{-1}  \P _{\a,j}^{-1} +   \P
 _{\a,i}^{-1} \P _{3,j}^{-1}\Big) \Big (\sum _{p=0}^n  \f_{,3}^p \cdot
 \f_{,\a}^{n-p}\Big) 
\end{equation}
and
\begin{equation}\label{p5} 
\dis\w  C_{ij}^n =\P _{3,i}^{-1}\P _{3,j}^{-1}\sum _{p=0}^n \f_{,3}^p
\cdot \f_{,3}^{n-p}, 
\end{equation}
where the derivatives of $\f^p$ are taken at point $z_{h}\circ \P
^{-1} (\w x)$. 
\end{lem}
Notice that the coefficients $\w{E}_{ij}^n$ depend on $h$.

\noindent{\bf Proof}
We have 
$$
 \w{E}_{ij}(h) = \frac{1}{2} (\w{\f}_{,i}\cdot \w{\f}_{,j} - \delta _{ij}).
$$
Writing the derivatives of $\w{\f}$ in term of those of $\f$ we get,
$$
 \w{\f}_{i,j} (\w x) = \Big(\f_{i,k}\circ z_{h} \circ \P^{-1}(\w
 x)\Big)\Big( (z_{h})_{k,l}\circ \P^{-1}(\w x)\Big) \P_{l,j}^{-1}(\w
 x), 
$$
with
$$
(z_{h})_{\alpha,\beta}=\delta_{\alpha\beta}\;\text{et}\;(z_{h})_{i,3}=
(z_{h})_{3,i}=\h\delta_{i3}. 
$$
Since $\f(x) = \sum_{p\geq 0}h^p \f^p (x)$, we obtain
$$
\w{\f}_{,j} (\w x) = \sum_{p\geq 0}h^p \Big(\Big( \big(\f_{,\a}^p\circ
z_{h} \circ \P^{-1}(\w x)\big) \P_{\a,j}^{-1}(\w x) + \h
\big(\f_{,3}^p\circ z_{h} \circ \P^{-1}(\w x)\big) \P_{3,j}^{-1}(\w
x)\Big)\Big). 
$$
Taking the scalar product, we obtain 
\begin{multline*}
\w{\f}_{,i} (\w x)\cdot \w{\f}_{,j} (\w x) = \sum_{n\geq 0} h^n
\sum_{p=0}^n\Big(\big(\f_{,\a}^p\circ z_{h} \circ \P^{-1}(\w x)\big)
\P_{\a,i}^{-1}(\w x)  + \h\big( \f_{,3}^p\circ z_{h} \circ \P^{-1}(\w
x)\big) \P_{3,i}^{-1}(\w x)\Big)\\ 
 \cdot \Big(\big(\f_{,\beta}^{n-p}\circ z_{h} \circ \P^{-1}(\w x)\big)
 \P_{\beta,j}^{-1}(\w x) + \h \big(\f_{,3}^{n-p}\circ z_{h} \circ
 \P^{-1}(\w x)\big) \P_{3,j}^{-1}(\w x)\Big), 
\end{multline*}
and thus
$$
\w{\f}_{,i} (\w x)\cdot \w{\f}_{,j} (\w x)  = \sum_{n\geq 0}h^n \Big(
\big(A^n _{ij} + \h B^n _{ij}+ \hc C^n _{ij}\big)\Big). 
$$
with $\w A^n _{ij} ,\,  \w B^n _{ij} $ and $\w C^n _{ij}$ defined as
above, which gives the result. 
\qed
\begin{lem}
The elastic energy term $I(h)(\f)$ admits an asymptotic expansion in
powers of $h$ 
$$
 I(h)(\f) = \sum_{n=-4}^{+\infty}  I^n_h (\f)h^n, 
$$
where
$$
 I^{-4}_h(\f) = \int_{\O}C^{ijkl} E_{ij}^{-2} E_{kl}^{-2} c_0 \,dx,
$$
$$
I^{-3}_h(\f) = \sum _{p=-4}^{-3}\int_{\O}C^{ijkl} E_{ij}^{p+2}E_{kl}^{-p-5} c_0 \,dx
+ \int_{\O}C^{ijkl} E_{ij}^{-2} E_{kl}^{-2} c_1 \,dx 
$$
and for every $n\geq -2$
\begin{multline*}
I^n_h(\f) = \sum _{p=-4}^n \int_{\O}C^{ijkl} E_{ij}^{p+2}E_{kl}^{n-p-2} c_0 \,dx +
\sum _{p=-4}^{n-1}\int_{\O}C^{ijkl} E_{ij}^{p+2} E_{kl}^{n-p-3} c_1 \,dx \\
+ \sum _{p=-4}^{n-2}\int_{\O}C^{ijkl} E_{ij}^{p+2} E_{kl}^{n-p-4} c_2 \,dx.
\end{multline*}
where
\begin{equation}\label{p1} 
 {E}_{ij}^{-2} = \frac{1}{2}   C_{ij}^0,\, {E}_{ij}^{-1} =
 \frac{1}{2}(  B_{ij}^0 +   C_{ij}^1),\,  
  {E}_{ij}^{0} = \frac{1}{2}(  A_{ij}^0 +   B_{ij}^1 +   C_{ij}^2 -
  \delta _{ij})\end{equation} 
and for every $n>0$
\begin{equation}\label{p2}
 {E}_{ij}^{n} = \frac{1}{2}(  A_{ij}^n +   B_{ij}^{n+1} +   C_{ij}^{n+2} ),
\end{equation}
with 
\begin{equation}\label{p3}
  A_{ij}^n =\P _{\a,i}^{-1}\big(\P (z^{-1}(x))\big) \P _{\beta,j}^{-1}
  \big(\P (z^{-1}(x))\big)\sum _{p=0}^n \f_{,\a}^p(x) \cdot
  \f_{,\beta}^{n-p}(x)   , 
\end{equation}
\begin{multline}\label{p4}
   B_{ij}^n = \Big(  \P _{3,i}^{-1}\big(\P (z^{-1}(x))\big)  \P
   _{\a,j}^{-1}\big(\P (z^{-1}(x))\big) +   \P _{\a,i}^{-1}\big(\P
   (z^{-1}(x))\big) \P _{3,j}^{-1}\big(\P (z^{-1}(x))\big)\Big)\\ \Big
   (\sum _{p=0}^n  \f_{,3}^p(x)  \cdot \f_{,\a}^{n-p}(x) \Big) 
\end{multline}
and
\begin{equation}\label{p5} 
\dis   C_{ij}^n =\P _{3,i}^{-1}\big(\P (z^{-1}(x))\big)\P
_{3,j}^{-1}\big(\P (z^{-1}(x))\big)\sum _{p=0}^n \f_{,3}^p (x) \cdot
\f_{,3}^{n-p}(x), 
\end{equation}
where $c_0$, $c_1$, $c_2$ are defined in Lemma \ref{l1}.
\end{lem}
\noindent{\bf Proof}
We have that
$$
 \w{I}^h (\w{\f}) = \int _{\w{\Omega}_h} C^{ijkl}\w{E}_{ij}^h
 (\w{\f}^h) \w{E}_{kl}^h (\w{\f}^h) dx.  
$$
Replacing $ \w{E}_{ij}^h (\w{\f}^h)$ and $\w{E}_{kl}^h (\w{\f}^h)$ by
their asymptotic expansion in powers of $h$ and applying the 
product of the two series, we obtain
$$ 
\w{I}^h (\w{\f}) = \sum _{n\geq -4} \Big( \Big(\sum _{p=-4}^n \int
_{\w{\Omega}_h} C^{ijkl}\w{E}_{ij}^{p+2} (\w{\f}^h)\w{E}_{kl}^{n-p-2}
(\w{\f}^h) dx\Big) \Big) h^n. 
$$
Next, by changing variables we obtain
$$
 I(\f) =  \sum _{n\geq -4}  \Big( \Big(\sum
 _{p=-4}^n \int _{\O} C^{ijkl}E_{ij}^{p+2}E_{kl}^{n-p-2}\; \det \grd
 \P(x_1, x_2, hx_3) dx\Big) \Big)h^n.  
$$
Then, replacing the Jacobian by its asymptotic expansion in powers of $h$,
we obtain
$$
 I(\f) =  \sum _{n\geq -4}  \Big( \Big(\sum _{p=-4}^n \int _{\O} C^{ijkl}
E_{ij}^{p+2}E_{kl}^{n-p-2} (c_0 +hc_1 +h^2 c_2) dx \Big) \Big)h^n,
$$
which gives the result.
\qed 
\begin{lem}\label{toy1}
The second order energy term $K(h)(\f)$ admits an asymptotic expansion
in powers of $h$ 
$$
 K(h)(\f)= \sum _{n=-4} ^{+\infty}  K^n_h (\f) h^n 
$$
with
\begin{align*}
 K^{-4}_h(\f) &= \int _{\O}   c_0\, Q(d_{ij}^{-2}\pt e_i \pt e_j) dx, \\
 K^{-3}_h(\f) &= \int _{\O} \big(  c_1\, Q(d_{ij}^{-2}\pt e_i \pt e_j) +
\sum _{p=-4}^{-3} c_0\, B(d_{ij}^{p+2} \pt e_i\pt e_j , d_{ij}^{-3-p-2}\pt e_i\pt e_j)\big)dx
\end{align*}
and for every $n\geq -2$
\begin{multline*}
K^n_h(\f)= \int_{\O}\Big( c_2 \sum _{p=-4}^{n-2} B(d_{ij}^{p+2} \pt e_i\pt e_j,d_{ij}^{n-p-4}\pt e_i\pt e_j)\\ +c_1 \sum _{p=-4}^{n-1} B(d_{ij}^{p+2}\pt e_i\pt e_j,
d_{ij}^{n-p-3}\pt e_i\pt e_j)
+ c_0 \sum _{p=-4}^n B(d_{ij}^{p+2}\pt e_i\pt e_j, d_{ij}^{n-p-2}\pt e_i\pt e_j)\Big) dx,  
\end{multline*}
where
\begin{equation}\label{d1}
 d_{ij}^{-2} = \z _{ij}^0,\;
d_{ij}^{-1} = \gamma _{ij}^0 +\z _{ij}^1 \,\text{ and for every } \;
n\geq 0, \;d_{ij}^n = \x _{ij}^n + \gamma _{ij}^{n+1} +\z _{ij}^{n+2}, 
\end{equation}
with
\begin{equation}\label{d2}
 \x _{ij}^{p}= \f_{,\a \beta}^p \P_{\a,i}^{-1} \P_{\beta,j}^{-1} +
 \f_{,\a}^p \P _{\a,ij}^{-1} , 
\end{equation}
\begin{equation}\label{d3}
 \gamma _{ij}^{p}=\f_{,\a3}^p(\P_{3,i}^{-1} \P_{\a,j}^{-1}  +
 \P_{\a,i}^{-1} \P_{3,j}^{-1} ) + \f_{,3}^p \P_{3,ij}^{-1} 
\end{equation}
and
\begin{equation}\label{d4}
 \z _{ij}^p = \f_{,33}^p \P_{3,i}^{-1} \P_{3,j}^{-1}.
\end{equation}

\end{lem}
In the previous formulas, the derivatives $\P _{k,i}^{-1}$ and $\P
_{k,ij}^{-1}$ are taken at $\P (x_1, x_2, hx_3)$ .
\noindent{\bf Proof}
We have 
$$
 \w{\f}(\w x) = \f\circ z_{h} \circ \P ^{-1}(\w x).
$$
Differentiating twice we obtain
\begin{multline*}
 \w{\f}_{,ij}(\w x) = \Big(\f_{,\a \beta}\circ z_{h}\circ \P ^{-1}(\w
 x)\Big) \P_{\a,i}^{-1}(\w x) \P_{\beta,j}^{-1}(\w x)  +\Big(
 \f_{,\a}\circ z_{h}\circ \P ^{-1}(\w x)\Big) \P_{\a,ij}^{-1}(\w x)
 \\ 
+  \h \Big[\Big( \f_{,\a 3}\circ z_{h}\circ \P ^{-1}(\w x)\Big)
\big(\P_{\a,i}^{-1}(\w x)  \P_{3,j}^{-1}(\w x) +  \P_{\a,j}^{-1}(\w x)
\P_{3,i}^{-1}(\w x)\big)+ \Big(\f_{,3}\circ z_{h}\circ \P ^{-1}(\w
x)\Big) \P_{3,ij}^{-1}(\w x)\Big]\\
 + \hc \Big(\f_{,33}\circ z_{h}\circ \P ^{-1}(\w x)\Big)
 \P_{3,i}^{-1}(\w x) \P_{3,j}^{-1}(\w x). 
\end{multline*}
Next, we replace $\f$ by its asymptotic expansion in powers of $h$, 
we obtain
$$
 \w{\f}_{,ij} = \sum _ {p\geq -2} \w{d}_{ij}^p h^p,
$$
with 
$$
 \w{d}_{ij}^{-2} = \w{\z} _{ij}^0,\; \w{d}_{ij}^{-1} = \w{\gamma}
 _{ij}^0 +\w{\z} _{ij}^1\;\text{et pour tout}\; n\geq 0,
 \;\w{d}_{ij}^n = \w{\x} _{ij}^n + \w{\gamma} _{ij}^{n+1} +\w{\z}
 _{ij}^{n+2}, 
$$
where $ \w{\x} _{ij}^{p}$, $ \w{\gamma} _{ij}^{p}$ and $ \w{\z}
_{ij}^p$ have the same algebraic expressions respectively as $ {\x}
_{ij}^{p}$, $ {\gamma} _{ij}^{p}$ and ${\z} _{ij}^p$ defined in
(\ref{d2}) (\ref{d3}) and (\ref{d4}) noticing that
$\f_{,i}$, $\f_{,ij}$ are taken at $z_{h}\circ \P ^{-1}(\w x)$ and $
\P _{k,i}^{-1}$, $\P _{k,ij}^{-1}$ are taken at $\w x$. 
Next, since we can write $\nabla^2 \w{\f} = \w{\f}_{,ij}\pt e_i\pt e_j $, we obtain that
$$
Q (\nabla^2 \w{\f}) = \sum_{n\geq -4}  \Big( \sum _{p=-4}^n
 B( \w{d}_{ij}^{p+2}\pt e_i\pt e_j , \w{d}_{ij}^{n-p-2} \pt e_i\pt e_j )\Big)h^n. 
$$
Finally, by changing the variables we obtain that
\begin{align*}
K(h)(\f) &= \h \w{K}^h(\w{\f}) = \h \int_{\w{\Omega}_h} Q (\nabla^2 \w{\f}) dx \\ 
&= \int _{\O}  \sum _{n\geq -4}  \Big( \Big( \sum _{p=-4}^n
B( d_{ij}^{p+2}\pt e_i\pt e_j , d_{ij}^{n-p-2} \pt e_i\pt e_j )\Big) (c_0 + hc_1 + h^2 c_2)
\Big)h^n\, dx, 
\end{align*}
which gives the result.
\qed 
{\bf Proof of Proposition \ref{prop1} }
It suffices to write that $J(h)(\f) = I(h)(\f) + K(h)(\f)$ and to apply
Lemmas 3.3 and 3.4, to obtain the result.
\qed 
\begin{rk}
The terms of the asymptotic expansion of $J(h)$ are not
independent on the thickness $h$ since the first and second derivatives of
$\dis \P^{-1}$ are taken at the point $\dis\p(x_1, x_2, hx_3)$. However, thanks to Remark \ref{quickli} we have the following
proposition.
\end{rk}

\begin{prop}\label{Propzah}
The energy $J(h)$ admits an asymptotic expansion of the form
\beq
J(h)(\f)=\sum_{n\geq-4}J^n(\f)h^n, 
\eeq
where the $J^n$ are independent on $h$.
\end{prop}
\noindent{\bf Proof}\\
We obtain this result using Proposition \ref{prop1} and Remark \ref{quickli}.
\qed
\begin{rk}
The expression of $J^n$ derives from that of $J^n_h$ using the asymptotic expansion of $\Psi^{-1}(\Psi(x_1,x_2,hx_3))$. First, it is obvious that $J^{-4}=J^{-4}_0$, where $J^n_0=J^n_h$ when $h=0$. Next, minimizing $J^{-4}$ we obtain that $J^{-4}=J^{-4}_h=0$. This implies that $J^{-3}=J^{-3}_0$. Similarly, we obtain that  $J^{-2}=J^{-2}_0$ minimizing  $J^{-3}$, also that $J^{-1}=J^{-1}_0$ minimizing  $J^{-2}$ and finally, that  $J^{0}=J^{0}_0$ minimizing  $J^{-1}$.
\end{rk}
\section{Solving the minimization problems}
We will use Proposition \ref{Pantz} to obtain the expression of $J^0$ following the mimimization of the energies $ J^{-4}\,,  J^{-3}\,, J^{-2}\,$ and $\; J^{-1}$. We have the following result.
\begin{prop}\label{css}
Minimizing the energies $ J^{-4}\,,  J^{-3}\,, J^{-2}\; \text{and}\; J^{-1}$, we obtain the following $J^0$ energy 
\begin{equation*}
J^0 (\f) = \int _{\O} \Big( \sum _{i,j=1}^3  Q( d_{ij}^0\pt e_i\pt e_j) +\sum _{i,j,k,l=1}^3  C^{ijkl} E_{ij}^0 E_{kl}^0 \Big)\, \det A(x) \, dx,
\end{equation*}
where $\f\in V_h$ verifies
\begin{equation}
\f^0_{,3}=0 \;\text{ and }\; \f^1_{,33}=0,
\end{equation} 
with
\begin{equation}
d_{ij}^0 = \f_{,\a \beta}^0\, \P_{\a,i}^{-1}\, \P_{\beta,j}^{-1} + 
\f_{,\a3}^1\,\big( \P_{\a,i}^{-1}\, \P_{3,j}^{-1} +\P_{3,i}^{-1}\, \P_{\a,j}^{-1}  \big)+\f_{,3}^1\, \P_{3,ij}^{-1} + \f_{,\a}^0\, \P_{\a,ij}^{-1}+\f_{,33}^2\, \P_{3,i}^{-1}\, \P_{3,j}^{-1},
\end{equation} 
\begin{equation}
E_{ij}^0 = \frac{1}{2}\Big[\big( \f_{,\a}^0\cdot \f_{,\beta}^0 \big)\,\P_{\a,i}^{-1}\, \P_{\beta,j}^{-1} + \big(\f_{,\a}^0\cdot \f_{,3}^1\big) \big(\P_{\a,i}^{-1}\, \P_{3,j}^{-1}  +\P_{3,i}^{-1}\, \P_{\a,j}^{-1} \big)+ \nor\f_{,3}^1\nor^2 \,\P_{3,i}^{-1}\, \P_{3,j}^{-1} -\,
\delta _{ij}\Big]
\end{equation}
and the functions $\P _{i,j}^{-1}$, $\P _{i,jk}^{-1}$ are taken at point $\P (x_1, x_2, 0)$.
\end{prop}

The proof of the Proposition is a consequence of the following lemmas and propositions.
\begin{lem}\label{leme-4}
The energy $J^{-4}$ is of the form
\begin{equation}
J^{-4}(\f) = \int_{\O}\big(Q( \f_{,33}^0\pt e_i\pt e_j )  {\P_{3,i}^{-1}}^2 {\P_{3,j}^{-1}}^2 + \frac{1}{4} C^{ijkl}  \nor \f_{,3}^0 \nor ^4\P_{3,i}^{-1} \P_{3,j}^{-1} \P_{3,k}^{-1} \P_{3,l}^{-1}\big) c_0 \;dx 
\end{equation}
where the functions $\P _{i,j}^{-1}$ are taken at point $\P (x_1, x_2, 0)$.
\end{lem}
\noindent{\bf Proof}\\
We use Proposition \ref{prop1} which implies that 
\begin{equation}\label{stef}
J_h^{-4}(\f) = \int_{\O} \Big( Q( d_{ij}^{-2}\pt e_i\pt e_j) + C^{ijkl} E_{ij}^{-2} E_{kl}^{-2}\Big) c_0 \,dx,
\end{equation}
where, using Lemma \ref{toy1}, we have 
\begin{equation}
d_{ij}^{-2}= \z _{ij}^0 = \f_{,33}^0 \P_{3,i}^{-1} \P_{3,j}^{-1}
\end{equation}
and using Lemma \ref{toy2}, we have
\begin{equation}
 E_{ij}^{-2}= \frac{1}{2} C_{ij}^0=\P _{3,i}^{-1}\P _{3,j}^{-1} \nor \f_{,3}^0 \nor^2,
\end{equation}
noticing that the functions $\P _{i,j}^{-1}$ are taken at point $\P (x_1, x_2, hx_3)$. Next, since we have 
\beq
h^4 J(h)(\f)=\sum_{n\geq 0}J^{n-4}_h(\f) h^n=\sum_{n\geq 0}J^{n-4}(\f) h^n,
\eeq
this implies that for every $\f\in Q_{-4}$ we have 
\beq
J^{-4}(\f)=  J_0^{-4}(\f),
\eeq
where $J^n_o=J^n_h$ for $h=0$, and thus the result.
\qed
Now that we have the expression of the energy $J^{-4}$, we can minimize it on $V_h$ following Pantz's method. We obtain.
\begin{prop}\label{trav}
Solving the first minimization problem we obtain  
\begin{center}
$Q_{-3} = \big\{ \f \in V_h ,\, \f_{,3}^0 = 0\big\}$.
\end{center}
Moreover,  $\dis J^{-4}(\f)=0$ on $Q_{-3}$.
\end{prop}
\noindent{\bf Proof}\\
Using Proposition \ref{Pantz} we have 
\begin{equation*}
  Q_{-3} =\big\{ \f \in V_h, J^{-4} (\f) = \inf _{\f'\in V_h} J^{-4}(\f') \big\}.
\end{equation*}
Next, using Lemma \ref{leme-4}, we notice that to minimize the energy $J^{-4}$ on $V_h$, it suffices to consider deformations $\f$ verifying  
\beq
\f^0_{,3}=0.
\eeq
For these deformations, we have 
\beq
J^{-4}(\f) =0,
\eeq
which gives the result since $J^{-4}(\f) \geq 0$ and $J^{-4}(\f) > 0$ if $\f^0_{,3}\neq 0.$
\qed
\begin{corol}
If $\f\in Q_{-3}$ then 
\beq
J^{-4}_h (\f) =0.
\eeq
\end{corol}
\noindent{\bf Proof}\\
The proof is obvious using (\ref{stef}). 
\qed
Next, we consider the second minimization problem. We have the following result.
\begin{prop}
The second minimization problem gives 
\begin{equation*}
 Q_{-2}= Q_{-3}.
\end{equation*}
\end{prop}
\noindent{\bf Proof}\\
Using Proposition \ref{Pantz}, we have 
\begin{equation*}
 Q_{-2}=\big\{ \f \in Q_{-3} , J^{-3} (\f) = \inf _{\f'\in Q_{-3}} J^{-3}(\f') \big\}
\end{equation*}
and using Proposition \ref{prop1}, we have
\begin{multline}\label{cro}
J^{-3}_h(\f) =  \sum _{p=-4}^{-3}\int_{\O} \Big( B(d_{ij}^{p+2}\pt e_i\pt e_j, d_{ij}^{-p-5}\pt e_i\pt e_j)+ C^{ijkl} E_{ij}^{p+2}E_{kl}^{-p-5} \Big)c_0 \\
+ \Big( Q( d_{ij}^{-2}\pt e_i\pt e_j) +C^{ijkl} E_{ij}^{-2} E_{kl}^{-2}\Big) c_1 \,dx. 
\end{multline}
Then, for $\f \in Q_{-3}$, Proposition \ref{trav} gives that $d_{ij}^{-2}=0$ and $ E_{ij}^{-2}=0$ which implies using (\ref{cro}) that 
\begin{equation}
J_h^{-3}(\f) = 0.
\end{equation}
Since we already proved that 
\begin{equation}
J_h^{-4}(\f) = J^{-4}(\f)= 0,
\end{equation}
we obtain 
\beq
h^3 J(h)(\f)=\sum_{n\geq 0} J^{n-3}_h(\f)h^n=\sum_{n\geq 0} J^{n-3}(\f)h^n.
\eeq
Thus, for every $\f\in Q_{-3}$ we have
\beq
J^{-3}(\f) =J_0^{-3}(\f) = 0,
\eeq
which gives the result.
\qed
Next, we solve the following minimization problem that will provide a condition on the second derivative of $\dis \f^1$.
\begin{lem}\label{cab}
For every $\dis \f\in Q_{-2}$ we have
\begin{equation}
 J^{-2}(\f) = \int _{\O}Q( \f_{,33}^1\pt e_i\pt e_j) {\P_{3,i}^{-1}}^2 {\P_{3,j}^{-1}}^2 \det (a_1\nor a_2\nor a_3)dx,
\end{equation}
where the functions $\P _{i,j}^{-1}$ are taken at point $\P (x_1, x_2, 0)$.
\end{lem}
\noindent{\bf Proof}\\

Using Proposition \ref{prop1}, we have
\begin{multline}
J^{-2}_h (\f) = \sum _{p=-4}^{-2} \int_{\O} \Big(B( d_{ij}^{p+2}\pt e_i\pt e_j, d_{ij}^{-p-4}\pt e_i\pt e_j) +C^{ijkl} E_{ij}^{p+2}E_{kl}^{-p-4}\Big) c_0 \,dx\\
 +\sum _{p=-4}^{-3}\int_{\O} \Big(B(d_{ij}^{p+2} \pt e_i\pt e_j, d_{ij}^{-p-5}\pt e_i\pt e_j) +C^{ijkl} E_{ij}^{p+2} E_{kl}^{-p-5}\Big) c_1 \,dx\\
  +\int_{\O}\Big( Q(d_{ij}^{-2}\pt e_i\pt e_j ) +C^{ijkl} \nor E_{ij}^{-2}\nor^2 \Big) c_2 \,dx.
\end{multline}
When $\dis \f \in Q_{-2}$, that is $\dis \f^0_{,3}=0$, we know that $\dis d_{ij}^{-2}=0$ and $\dis E_{ij}^{-2}=0$. Thus, we have
\begin{equation}\label{ess}
J^{-2}_h(\f) = \int _{\O}\sum _{i,j,k,l=1}^3\big(Q(d_{ij}^{-1}\pt e_i\pt e_j) + C ^{ijkl} E_{ij}^{-1} E_{kl}^{-1}\big) \, \det (a_1\nor a_2\nor a_3)dx.
\end{equation}
Moreover, Lemma \ref{toy2} implies that
$$\begin{aligned}
 E_{ij}^{-1}= \frac{1}{2}( B_{ij}^{0} + C_{ij}^{1} )&=\frac{1}{2}\Big(\big(  \P _{3,i}^{-1}  \P _{\a,j}^{-1} +   \P _{\a,i}^{-1} \P _{3,j}^{-1}\big)  \f_{,3}^0 \cdot \f_{,\a}^0 + \P _{3,i}^{-1}\P _{3,j}^{-1}\sum _{p=0}^1 \f_{,3}^p \cdot \f_{,3}^{1-p} \Big)\\
&=\frac{1}{2}\Big(\big(  \P _{3,i}^{-1}  \P _{\a,j}^{-1} +   \P _{\a,i}^{-1} \P _{3,j}^{-1}\big)  \f_{,3}^0 \cdot \f_{,\a}^0 + \P _{3,i}^{-1}\P _{3,j}^{-1}\nor \f_{,3}^0 \cdot \f_{,3}^1\nor^2 \Big)
\end{aligned}$$
and Lemma \ref{toy1} implies that 
$$\begin{aligned}
d_{ij}^{-1} = \gamma _{ij}^0 +\z _{ij}^1=\f_{,\a3}^0(\P_{3,i}^{-1} \P_{\a,j}^{-1}  + \P_{\a,i}^{-1} \P_{3,j}^{-1} ) + \f_{,3}^0 \P_{3,ij}^{-1}+ \f_{,33}^1 \P_{3,i}^{-1} \P_{3,j}^{-1}.
\end{aligned}$$
Thus, when $\dis \f \in Q_{-2}$ we have 
\begin{center}
$\dis E_{ij}^{-1}=0$ and $\dis d_{ij}^{-1} =\f_{,33}^1 \P_{3,i}^{-1} \P_{3,j}^{-1}$,
\end{center}
where the functions $\P _{i,j}^{-1}$ are taken at point $\P (x_1, x_2, hx_3)$. Replacing these expressions into (\ref{ess}), we obtain that
\begin{equation}\label{anna}
 J^{-2}_h(\f) = \int _{\O}Q( \f_{,33}^1\pt e_i\pt e_j ) {\P_{3,i}^{-1}}^2 {\P_{3,j}^{-1}}^2 \det (a_1\nor a_2\nor a_3)dx,
\end{equation}
where the functions $\P _{i,j}^{-1}$ are taken at point $\P (x_1, x_2, hx_3)$. Next, since $\f\in Q_{-2}$ we have
\beq
h^2 J(h)(\f)=\sum_{n\geq 0} J^{n-2}_h(\f)h^n=\sum_{n\geq 0} J^{n-2}(\f)h^n,
\eeq
and thus, for every $\f\in Q_{-2}$ we have
\beq
J^{-2}(\f) =J_0^{-2}(\f),
\eeq
which gives the result.
\qed
As a consequence of Lemma \ref{cab}, we obtain the following Proposition.
\begin{prop}
The third minimization problem gives 
\begin{equation}
Q_{-1} = \big\{ \f \in V_h ,\, \f_{,3}^0 = 0\,\text{ and } \, \f_{,33}^1 = 0\big\}.
\end{equation}
\end{prop}
\noindent{\bf Proof}\\
Using Proposition \ref{Pantz} we have
\begin{equation*}
  Q_{-1} =\big\{ \f \in Q_{-2}, J^{-2} (\f) = \inf _{\f'\in Q_{-2} } J^{-2} (\f') \big\}.
\end{equation*}
Next, using Lemma \ref{cab} and since the quadratic form is positive, we know that in order to minimize the energy $\dis J^{-2}$, it is necessary and sufficient to consider $ \f$ such that $ \dis \f_{,33}^1=0$, which gives the result.
\qed
\begin{corol}
For every $\f\in  Q_{-1}$ we have 
\beq
 J^{-2}_h (\f)=0.
\eeq
\end{corol}
\noindent{\bf Proof}\\
We obtain the result using (\ref{anna}).
\qed
Next, we pass to the following minimization problem. We have the following Proposition.
\begin{prop}
The fourth minimization problem is trivial and gives
\begin{equation*}
Q_0= Q_{-1}.
\end{equation*}
\end{prop}
\noindent{\bf Proof}\\
Using Proposition \ref{Pantz} we have  
\begin{equation*}
  Q_0 =\big\{ \f \in Q_{-1}, J^{-1} (\f) = \inf _{\f'\in Q_{-1} } J^{-1} (\f') \big\}.
\end{equation*}
Moreover, using Proposition \ref{prop1}, we have
\begin{multline}
J^{-1}_h (\f) = \sum _{p=-4}^{-1} \int_{\O} \Big(B(d_{ij}^{p+2}\pt e_i\pt e_j, d_{ij}^{-p-3}\pt e_i\pt e_j) +C^{ijkl} E_{ij}^{p+2}E_{kl}^{-p-3}\Big) c_0 \,dx \\ +
\sum _{p=-4}^{-2}\int_{\O} \Big(B(d_{ij}^{p+2}\pt e_i\pt e_j, d_{ij}^{-p-4}\pt e_i\pt e_j) 
+C^{ijkl} E_{ij}^{p+2} E_{kl}^{-p-4}\Big) c_1 \,dx \\+ \sum _{p=-4}^{-3}\int_{\O}\Big( B(d_{ij}^{p+2}\pt e_i\pt e_j, d_{ij}^{-p-5}\pt e_i\pt e_j) +C^{ijkl} E_{ij}^{p+2} E_{kl}^{-p-5}\Big) c_2 \,dx.
\end{multline}
Next, using (\ref{d1}) and (\ref{p1}), we have that for every $\f \in  Q_{-1}$, 
\beq
d_{ij}^{-2}=d_{ij}^{-1}=0\;\text{ and}\; E_{ij}^{-2}=E_{ij}^{-1}=0,
\eeq
which implies that for every $ \f \in Q_{-1}$ we have 
\begin{equation}
J^{-1}_h (\f) =0.
\end{equation}
Since, for every $ \f \in Q_{-1}$ we have
\begin{equation}
J^{-2} (\f) = J^{-2}_h (\f) =0,
\end{equation}
it implies that for every $\f\in Q_{-1}$ we have
\beq
J^{-1}(\f)=J^{-1}_0 (\f)=0,
\eeq
which gives the result.
\qed
Finally, the following minimization problem will provide the proof of Proposition \ref{css}.\\
{\bf Proof of Proposition \ref{css}}\\

Using Proposition \ref{prop1} we have
\begin{multline}
J^0_h (\f) = \sum _{p=-4}^0 \int_{\O} \Big(B( d_{ij}^{p+2}\pt e_i\pt e_j, d_{ij}^{-p-2}\pt e_i\pt e_j) +C^{ijkl} E_{ij}^{p+2}E_{kl}^{-p-2}\Big) c_0 \,dx \\
+\sum _{p=-4}^{-1}\int_{\O} \Big(B(d_{ij}^{p+2}\pt e_i\pt e_j, d_{ij}^{-p-3}\pt e_i\pt e_j) 
+C^{ijkl} E_{ij}^{p+2} E_{kl}^{-p-3}\Big) c_1 \,dx \\
+ \sum _{p=-4}^{-2}\int_{\O}\Big( B(d_{ij}^{p+2}\pt e_i\pt e_j, d_{ij}^{-p-4}\pt e_i\pt e_j) +C^{ijkl} E_{ij}^{p+2} E_{kl}^{-p-4}\Big) c_2 \,dx.
\end{multline}
Cancelling the terms $d_{ij}^{-2},\,d_{ij}^{-1},\,E_{ij}^{-2}$ and $E_{ij}^{-1}$ we obtain that
\begin{equation}
J^0_h (\f) =  \int_{\O} \Big(Q( d_{ij}^0\pt e_i\pt e_j)  +C^{ijkl} E_{ij}^0 E_{kl}^0\Big) c_0 \,dx .
\end{equation}
Next, using Lemma \ref{toy2}, we have
\begin{multline}
 E_{ij}^0= \frac{1}{2}(A_{ij}^0 + B_{ij}^{1} + C_{ij}^{2}-\delta_{ij} )\\
=\frac{1}{2}\Big(\P _{\a,i}^{-1} \P _{\beta,j}^{-1}  \f_{,\a}^0 \cdot \f_{,\beta}^0+ \Big(  \P _{3,i}^{-1}  \P _{\a,j}^{-1} +   \P _{\a,i}^{-1} \P _{3,j}^{-1}\Big) \Big (\sum _{p=0}^1  \f_{,3}^p \cdot \f_{,\a}^{1-p}\Big) +\P _{3,i}^{-1}\P _{3,j}^{-1}\sum _{p=0}^2 \f_{,3}^p \cdot \f_{,3}^{2-p}-\delta_{ij}\Big)
\end{multline}
which gives for $\dis \f\in Q_0$
\begin{multline}
 E_{ij}^0= \frac{1}{2}\Big(\P _{\a,i}^{-1} \P _{\beta,j}^{-1}  \f_{,\a}^0 \cdot \f_{,\beta}^0+ \Big(  \P _{3,i}^{-1}  \P _{\a,j}^{-1} +   \P _{\a,i}^{-1} \P _{3,j}^{-1}\Big)   \f_{,3}^1 \cdot \f_{,\a}^0+\P _{3,i}^{-1}\P _{3,j}^{-1}\nor \f_{,3}^1 \nor^2-\delta_{ij} \Big).
\end{multline}
Then, using Lemma \ref{toy1}, we have
\begin{multline}
 d_{ij}^0= \x _{ij}^0 + \gamma _{ij}^{1} +\z _{ij}^{2}\\
=\f_{,\a \beta}^0 \P_{\a,i}^{-1} \P_{\beta,j}^{-1} + \f_{,\a}^0 \P _{\a,ij}^{-1}+\f_{,\a3}^1(\P_{3,i}^{-1} \P_{\a,j}^{-1}  + \P_{\a,i}^{-1} \P_{3,j}^{-1} ) + \f_{,3}^1 \P_{3,ij}^{-1}+\f_{,33}^2 \P_{3,i}^{-1} \P_{3,j}^{-1},
\end{multline}
where the functions $\P _{i,j}^{-1}$ and $\P _{i,jk}^{-1}$ are taken at point $\P (x_1, x_2,hx_3)$.
Since for $\dis \f\in Q_0$ we have
\beq
J^{-1}_h(\f)=J^{-1}(\f)=0,
\eeq
it implies that
\beq
J(h)(\f)=\sum_{n\geq0} J^n_h(\f)h^n=\sum_{n\geq0} J^n(\f)h^n.
\eeq
Then, we obtain that
\beq
J^{0}(\f) =J_0^{0}(\f),
\eeq
which completes the proof of Proposition \ref{css}.
\qed
\begin{rk}
The limit model obtained is in agreement with the one obtained in \cite{ledzor,ledzor2} using $\Gamma$-convergence techniques for the case of curved martensitic thin films where the quadratic form is the square of the norm.
\end{rk}
\begin{rk}
Applying our results for the case of planar films i.e. for $\p(x)=x
$ we recover the formal asymptotic expansion of
Pantz for the elastic energy of a  Saint Venant-Kirchhoff material and the limit interfacial energy part obtained by Bhattacharya and James for the martensitic thin films again when the quadratic form is the square of the norm.
\end{rk}
\section{Application to the planar, cylindrical and spherical cases}
In this section we will illustrate our results for films that are planar, cylindrical or spherical. In the cylindrical and spherical case, we will suppose for simplicity that the quadratic form is the square of the norm in $\mathbb{M}_{333}$.
\subsection{Case of planar films }
In this section we will apply our results for the case of planar films, i.e. $\P(x) = x $. The expression of the energy $e^0$ becomes
\begin{multline}
J^0(\f) = \int_{\O}\Big \{ \frac{\lambda}{4} \Big ( \nor \f^0_{,1}\nor ^2 + \nor \f^0_{,2}\nor ^2 +\nor \f^1_{,3}\nor ^2 - 3 \Big )^2 \\
+ \frac{\mu}{2} \Big [ (\nor \f^0_{,1}\nor ^2 -1)^2 +  (\nor \f^0_{,2}\nor ^2 -1)^2 +  (\nor \f^1_{,3}\nor ^2 -1)^2 
+2 (\f^0_{,1} \cdot \f^0_{,2} )^2 + 2 (\f^0_{,1} \cdot \f^1_{,3} )^2 
+2 (\f^0_{,2} \cdot \f^1_{,3} )^2 \Big ]  \\
+  \Big (Q( \f^0_{,\alpha\beta} \pt e_\alpha\pt e_\beta ) + Q( \f^1_{,\alpha 3} \pt e_\alpha \pt e_3 )+ Q( \f^1_{,\alpha 3} \pt e_3 \pt e_\alpha ) +  Q( \f^2_{,3 3} \pt e_3 \pt e_3 )   \Big ) \Big \}  dx.
\end{multline} 
\begin{rk}
In the case of planar films, the derivatives of $\P$ are constants and thus, the dependence on $h$ that appears in the curved case following the first five minimization problems does not exist. Consequently, we can go further in solving the minimization problems.  
\end{rk}
We have the following result.
\begin{theorem}\label{te}
	Suppose that $Q_{i33}>0$ for $i=1,2,3$. Thus, minimizing $J^0$ on $\dis Q_0$ we obtain that 
\begin{multline}
Q_1=\big\{ \f\in V_h\;\text{verifying}\; \f^0_{,3}=0,\; \f^1_{,33}=0,\; J^0_1( \f^0,\f^1)=\inf_{(\xi,\eta) \in Q^0_2} J^0_1( \xi,\eta), \; \f^2_{,33}=0 \\
\;\text{and}\;\f^0 (x) = (\w{a_1}\nor\w{a_2}\nor 0) x ,\, \f^1_{,3}(x)=x_3 \w{a_3} \;\text{on} \; \partial \o \times (-\frac{1}{2},\frac{1}{2}) \big\},
\end{multline}
with $\w A=(\w a_1\nor \w a_2\nor \w a_3)$ being the constant matrix appearing in the boundary condition, 
\begin{multline}
 Q^0_2=  \Big\{ (z_0,z_1)\in  H^2 (\O ;  \mathbb{R}^3)\times H^2 (\O ;  \mathbb{R}^3) \; \text{such that there exist}\; (\f^2, \f^3,...) \; \text{verifying}\;\\
 z_0+ h z_1+ h^2 \f^2 +...\in Q^0\Big \}
\end{multline}
and
\begin{multline}\label{dans}
J^0_1(\f^0, \f^1) = \int_{\O}\Big \{  \Big (Q( \f^0_{,\alpha\beta} \pt e_\alpha\pt e_\beta ) + Q( \f^1_{,\alpha 3} \pt e_\alpha \pt e_3 )+ Q( \f^1_{,\alpha 3} \pt e_3 \pt e_\alpha )    \Big )\\
	+\frac{\lambda}{4} \Big ( \nor \f^0_{,1}\nor ^2 + \nor \f^0_{,2}\nor ^2 +\nor \f^1_{,3}\nor ^2 - 3 \Big )^2 \\
	+ \frac{\mu}{2} \Big [ (\nor \f^0_{,1}\nor ^2 -1)^2 +  (\nor \f^0_{,2}\nor ^2 -1)^2 +  (\nor \f^1_{,3}\nor ^2 -1)^2 
	+2 (\f^0_{,1} \cdot \f^0_{,2} )^2 + 2 (\f^0_{,1} \cdot \f^1_{,3} )^2 
	+2 (\f^0_{,2} \cdot \f^1_{,3} )^2 \Big ]\Big \}  dx.
\end{multline}
\end{theorem}

\noindent{\bf Proof}\\

Recall that 
\begin{multline}
Q^0 = \Big \{ \f = \sum_{i\geq 0}h^i\f^i\;\text{ with}  \; \f^n \in  H^2 (\O ;  \mathbb{R}^3), \; \; \f^0_{,3} =0, \; \; \f^1_{,33} =0 \; 
\\ \text{and}  \; \f (x) = (\w{a_1}\nor\w{a_2}\nor h \w{a_3}) x \;\text{on} \; \partial \o \times (-\frac{1}{2},\frac{1}{2})\Big\},    
\end{multline}
where $\dis \w{A} = (\w{a_1}\nor\w{a_2}\nor \w{a_3})$ is a constant matrix in $\M_{33}$. The boundary condition on $\dis \partial \o\times (-\frac{1}{2},\frac{1}{2})$ writes
\begin{equation}\label{bord}
 \f^0 (x) = ( \w{a_1} \nor \w{a_2} \nor 0) x,\, \f^1 (x) = (0 \nor 0 \nor \w{a_3}) x = x_3 \w {a_3}\;\text{and}\;  \f^n (x) =0 \;\text{for}\; n \geq 2.
\end{equation}
Let
\begin{equation}
J^0_1(\f^0, \f^1) = \inf_{ z^2 \in Q^0_1} J^0(\f^0, \f^1, z^2), 
\end{equation}
with 
\begin{equation}
 Q^0_1 = \Big\{ z\in  H^2 (\O ;  \mathbb{R}^3) \; \text{such that there exist}\; (\f^0, \f^1, \f^3, ...) \; \text{verifying}\; \f^0+ h \f^1+ h^2 z +h^3 \f^3+...\in Q^0\Big\}.
\end{equation}
We have that for every $\dis \f^2$ verifying $\displaystyle  \f^2_{,33} = 0$,
we obtain (\ref{dans} ).
\qed
\begin{rk}
Notice that the result is in accordance with the limit model obtained by Bhattacharya and James (voir\cite{bhajam}) for a martensitic plate obtained by $\Gamma$ convergence techniques. 
\end{rk}
\begin{rk}
	The hypothesis $q_{i33}>0 $ for $i=1,2,3$ is not necessary to obtain (\ref{dans}). Indeed, if one of the $q_{i33}$ is equal to zero, then the corresponding term do not appear in the energy from the beginning. 
\end{rk}
Finally, in the case when the quadratic form represents the square of the norm in $\mathbb{M}_{333}$, let us define the two-dimensional energy $J$ by
\begin{multline}
J(u) = \int_{\o}\Big \{ 2 \nor \grd u\nor^2 +\frac{\lambda}{4} \Big ( \nor \f^0_{,1}\nor ^2 + \nor \f^0_{,2}\nor ^2 +\nor u\nor ^2 - 3 \Big )^2 \\
+ \frac{\mu}{2} \Big [ (\nor \f^0_{,1}\nor ^2 -1)^2 +  (\nor \f^0_{,2}\nor ^2 -1)^2 +  (\nor u\nor ^2 -1)^2 
+2 (\f^0_{,1} \cdot \f^0_{,2} )^2 + 2 (\f^0_{,1} \cdot u )^2 +2 (\f^0_{,2} \cdot u )^2 \Big ]\Big \}  dx.
\end{multline}
Notice that the minimizer of the energy $J$ is the $\f^1_{,3}$ that minimizes the energy $J^0$. Writing the corresponding Euler-Lagrange equation we obtain that the minimizers are solutions of the following partial differential equation
\begin{displaymath}
\begin{cases}
[\lambda (\nor \f^0_{,1} \nor ^2+  \nor \f^0_{,2} \nor ^2 +\nor u \nor^2 -3)+ 2 \mu(\nor u \nor ^2-1)]u+2\mu[(u,\f^0_{,1})\f^0_{,1}+(u,\f^0_{,2})\f^0_{,2}]\\
-4 \Delta u=0\; \text{in}\; \o\\
u=\w A a_3\; \text{on} \; \partial \o.
\end{cases}
\end{displaymath} 
\subsection{Case of a cylindrical film}
In order to illustrate our results on an explicit example of curved film, we apply them in the case of a thin film of cylindrical
form, length $1$ and of radius $r$ when the quadratic form represents the square of the norm. We recall that in the case of a
portion of a cylinder, the diffeomorphism $\p$ is of the form:
\begin{eqnarray*}
\p :\omega := ]0;1[ \times ]0;\frac{\pi}{2}[ &\cv& \mathbb{R}^3\\
(x_1, x_2)&\mapsto& (x_1,r \cos x_2, r \sin x_2 ),
\end{eqnarray*}
which implies that 
\begin{center}
$a_1(x_1, x_2) = (1, 0 ,0)$, $a_2(x_1, x_2) = (0, - r\sin x_2, r\cos x_2)$ and $a_3(x_1, x_2) = (0, - \cos x_2,-\sin x_2) $.
\end{center}
Thus, $c_0 =r$ and $\P$ writes
\begin{eqnarray*}
\P : \Omega_1:= ]0;1[ \times ]0;\pi[\times ]-\frac{h}{2};\frac{h}{2} [ &\cv& \widetilde{\Omega}_{h}\\
(x_1, x_2, x_3)&\mapsto& (x_1, (r- x_3)\cos x_2,(r-x_3) \sin x_2 ),
\end{eqnarray*}
which gives that 
\begin{eqnarray*}
\P^{-1} :  \widetilde{\Omega}_{h} &\cv& \O \\
(z_1, z_2, z_3)&\mapsto& \Big(z_1,\arctan \frac {z_3}{z_2},r-\sqrt{z_2 ^2 + z_3 ^2 } \Big).
\end{eqnarray*}
Next, we compute $\dis \P^{-1}$ and its derivatives at the point 
\begin{center}
$\dis \P (x_1, x_2, 0)=(x_1,r\cos x_2, r \sin x_2) $.
\end{center}
We obtain the values of $d_{ij}^0$ and $E_{ij}^0$, which give the following energy $J^0$
\begin{multline} 
J^0(\f^0, \f^1, \f^2) = \int_{\O}\Big \{\frac{\lambda}{4} \Big ( \nor \f^0_{,1}\nor ^2 + \frac{1}{r}\nor \f^0_{,2}\nor ^2 +\nor \f^1_{,3}\nor ^2 - 3 \Big )^2\\
+ \frac{\mu}{2} \Big [ (\nor \f^0_{,1}\nor ^2 -1)^2 +  (\frac{1}{r}\nor \f^0_{,2}\nor ^2 -1)^2 +  (\nor \f^1_{,3}\nor ^2 -1)^2 +2 (\f^0_{,1} \cdot \f^0_{,2} )^2 + 2 (\f^0_{,1} \cdot \f^1_{,3} )^2 +2 (\f^0_{,2} \cdot \f^1_{,3} )^2 \Big ]\\
+  \Big (\nor \f^0_{,11}\nor^2 +\frac{2}{r^2} \nor \f^0_{,12}\nor^2 + 2  \nor\f^1_{,13}\nor^2+ \big\nor\frac{1}{r^2}\f^0_{,22}-\frac{1}{r}\f^1_{,3}\big\nor^2 +2 \nor \frac{1}{r}\f^1_{,23} +\frac{1}{r^2}\f^0_{,2} \nor^2 + \nor\f^2_{,33}\nor^2 \Big ) \Big \}\, r dx.
\end{multline}
Notice that we can minimize $J^0$ with respect to $\dis \f^2_{,33}$ taking
\beq
\dis \f^2_{,33} = 0.
\eeq
The energy becomes
\begin{multline} 
J^0(\f^0, \f^1) = \int_{\O}\Big \{\frac{\lambda}{4} \Big ( \nor \f^0_{,1}\nor ^2 + \frac{1}{r}\nor \f^0_{,2}\nor ^2 +\nor \f^1_{,3}\nor ^2 - 3 \Big )^2\\
+ \frac{\mu}{2} \Big [ (\nor \f^0_{,1}\nor ^2 -1)^2 +  (\frac{1}{r}\nor \f^0_{,2}\nor ^2 -1)^2 +  (\nor \f^1_{,3}\nor ^2 -1)^2 +2 (\f^0_{,1} \cdot \f^0_{,2} )^2 + 2 (\f^0_{,1} \cdot \f^1_{,3} )^2 +2 (\f^0_{,2} \cdot \f^1_{,3} )^2 \Big ]\\
+  \Big (\nor \f^0_{,11}\nor^2 +\frac{2}{r^2} \nor \f^0_{,12}\nor^2 + 2  \nor\f^1_{,13}\nor^2+ \big\nor\frac{1}{r^2}\f^0_{,22}-\frac{1}{r}\f^1_{,3}\big\nor^2 +2 \nor \frac{1}{r}\f^1_{,23} +\frac{1}{r^2}\f^0_{,2} \nor^2  \Big ) \Big \}\, r dx.
\end{multline}
The main difference, at this stage, with the case of plates,
appears in the existence of terms of order $1$ in the second order energy part, namely $\frac{1}{r^2}\f^0_{, 2 } $ and
$\frac{1}{r}\f^1_{, 3}$.  This makes the minimization of $J^0$ with respect to $\dis \f^1$ clearly different from the case of planar films.  Indeed, in order to minimize our energy with respect to $ \f^1$,
we consider the energy defined on $H^1(\o;\R^3)$ with $ \f^0_{, 1}$, $ \f^0_{, 2}$ and $\f^0_{,
  22}$ constant, writing
\begin{multline} 
 J(u) = \int_{\o}\Big \{\frac{\lambda}{4} \Big ( \nor {\f^0_{,1}} \nor ^2 + \frac{1}{r}\nor {\f^0_{,2}} \nor ^2 +\nor u \nor ^2 - 3 \Big )^2 + \frac{\mu}{2} \Big [ (\nor u \nor ^2 -1)^2 + 2 ({\f^0_{,1}} \cdot u )^2 +2 ({\f^0_{,2}} \cdot u )^2 \Big ]\\
+ \Big (2 \nor u_{,1}\nor^2+ \nor \frac{1}{r} u -\frac{1}{r^2} \f^0_{,22} \nor^2 + 2 \nor \frac{1}{r}u_{,2}+\frac{1}{r^2}{\f^0_{,2}}\nor^2\Big ) \Big \} \, r dx.
\end{multline}  
Notice that the minimizer of the last energy is the same $\f^1_{,3}$ that minimizes $J^0$. Considering the corresponding Euler-Lagrange equation we obtain that the minimizers verifies the following partial differential equations
\begin{displaymath}
\begin{cases}
[\lambda (\nor \f^0_{,1} \nor ^2+ \r \nor \f^0_{,2} \nor ^2 +\nor u \nor^2 -3)+ 2 \mu(\nor u \nor ^2-1)]u+2\mu[(u,\f^0_{,1})\f^0_{,1}+(u,\f^0_{,2})\f^0_{,2}]\\
+2\frac{1}{r} \big (\frac{1}{r} u -\frac{1}{r^2} \f^0_{,22} \big)-4 \w\Delta u=0\; \text{in}\; \o\\
u=\w A a_3\; \text{on} \; \partial \o,
\end{cases}
\end{displaymath} 
where
\beq
 \w\Delta u= u_{,11}+\frac{1}{r^2}u_{,22}.
\eeq

\subsection{Case of a portion of a sphere}
We consider the case of a thin film of spherical form, which is
rather different from that of the cylinder and further away from the case of planar films. Supposing that the quadratic form is the square of the norm, let us consider a spherical midsurface for which the diffeomorphism $\p
$ reads :
\begin{eqnarray*}
\p :\omega :=I_1\times I_2 \subset \mathbb{R}^2  &\cv& \widetilde{S}\subset \mathbb{R}^3\\
(x_1, x_2)&\mapsto& (\sin x_1 \cos x_2, \sin x_1 \sin x_2, \cos x_1 )
\end{eqnarray*}
where $x_1$ and $x_2$ represents the latitudinal and longitudinal coordinates, $I_1$ and $I_2$ are open intervals of $\mathbb{R}$ such that $\widetilde{S}$ do not contain the poles of the sphere in order to avoid singularities. We obtain that 
\begin{align*}
a_1(x_1, x_2) &= (\cos x_1 \cos x_2, \cos x_1 \sin x_2 , -\sin x_1),\\
a_2(x_1, x_2 )&= ( -\sin x_1 \sin x_2, \sin x_1 \cos x_2 , 0)\\
\text{and } a_3(x_1, x_2) &= (\sin x_1 \cos x_2, \sin x_1 \sin x_2,\cos x_1) .
\end{align*}
Thus, $c_0 = \sin x_1$ and $\P$ reads
\begin{eqnarray*}
\P : \Omega_1:=I_1 \times I_2 \times ]-\frac{h}{2};\frac{h}{2} [ &\cv& \widetilde{\Omega}_{h}\\
(x_1, x_2, x_3)&\mapsto& \big((1+x_3)\sin x_1 \cos x_2, (1+x_3)\sin x_1 \sin x_2, (1+x_3)\cos x_1 \big) ,
\end{eqnarray*}
which implies that 
\begin{eqnarray*}
\P^{-1} :  \widetilde{\Omega}_{h} &\cv& \O \\
(z_1, z_2, z_3)&\mapsto& \Big(\arctan \sqrt{\frac {z_1^2 + z_2^2}{z_3^2}},\arctan \frac {z_2}{z_1}, -1+\sqrt{z_1^2 +z_2 ^2 + z_3 ^2 } \;\Big).
\end{eqnarray*}
Then, we compute $\dis \P^{-1}$ and its derivatives at point $\dis \P (x_1, x_2, 0) $ to obtain the values of $d_{ij}^0$ and $\dis E_{ij}^0$, which gives the following energy $\dis J^0$,
$$\begin{aligned} 
&J^0(\f) = \int_{\O} \Big \{ \Big [ \nor \f^0_{,11}\nor ^2 +\frac{1}{\sin ^4 x_1}\nor \f^0_{,22}\nor ^2 + \frac{2}{\sin ^2 x_1}\nor \f^0_{,12}\nor ^2 +2 \nor \f^1_{,13}\nor ^2 + \frac{2}{\sin ^2 x_1} \nor \f^1_{,23}\nor ^2  + \nor \f^2_{,33}\nor ^2 \\
&+ \frac {2\cos x_1} {\sin^3 x_1} \; \f^0_{,22}   \cdot \f^0_{,1} -  \frac {4 \cos x_1}{\sin^3 x_1} \; \f^0_{,12}  \cdot \f^0_{,2} - 4 \; \f^1_{,13} \cdot \f^0_{,1} -   \frac {4 }{\sin^2 x_1} \;\f^1_{,23}  \cdot \f^0_{,2} \\
&+   \frac {2}{\sin^2 x_1} \; \f^0_{,22}  \cdot \f^1_{,3} + \nor \f^1_{,3}\nor ^2 +\nor \f^0_{,1}\nor ^2 \big( 4 \cos ^2 x_1 \sin ^2 x_1 + \cot ^2 x_1 +2 \big) + \nor \f^0_{,2}\nor ^2   \Big ] \\
&+ \frac{\mu}{2} \Big [ \big ( \nor \f^0_{,1}\nor^2 - 1  \big )^2 + \big (  \frac{1}{\sin ^2 x_1}\nor \f ^0_{,2}\nor^2 -1 \big)^2+ \big(  \nor \f^1_{,3}\nor^2 - 1 \big)^2 + \frac{2}{\sin^2 x_1} \big( \f^0_{,1}\cdot \f^0_{,2}\big )^2 + 2 \big( \f^0_{,1}\cdot \f^1_{,3}\big)^2 \\
&+  \frac{2}{\sin^2 x_1} \big( \f^0_{,2}\cdot \f^1_{,3}\big)^2 \Big ] + \frac{\lambda}{4} \Big(\nor \f^0_{,1}\nor^2 + \frac{1}{\sin ^2 x_1}\nor \f ^0_{,2}\nor^2 +\nor \f^1_{,3}\nor ^2  - 3  \Big )^2 \Big \} \sin x_1 dx.
\end{aligned}$$
We notice that here also we can minimize our energy with respect to $\f^2_{, 33}$ only by canceling it, which means that $\f^2 $ is linear with respect to the third variable.  Consequently, all the remaining terms
are independent of the third variable. This enables us to integrate
with respect to the third variable obtaining the following expression
$$\begin{aligned} 
&J^0(\f) = \int_{\o}  \Big \{  \Big [ \nor \f^0_{,11}\nor ^2  + \frac{2}{\sin ^2 x_1}\nor \f^0_{,12}\nor ^2 +2 \nor \f^1_{,13}-\f^0_{,1}\nor ^2 + \frac{2}{\sin ^2 x_1} \nor \f^1_{,23}-\f^0_{,2}\nor ^2  \\
&+ \frac {2\cos x_1} {\sin^3 x_1} \; \f^0_{,22}   \cdot \f^0_{,1} -  \frac {4 \cos x_1}{\sin^3 x_1} \; \f^0_{,12}  \cdot \f^0_{,2}  +   \nor\frac {1}{\sin^2 x_1} \; \f^0_{,22}+ \f^1_{,3}\nor^2 \\
& +\nor \f^0_{,1}\nor ^2 \big( 4 \cos ^2 x_1 \sin ^2 x_1 + \cot ^2 x_1 +1 \big) + \nor \f^0_{,2}\nor ^2(1- \frac{2}{\sin ^2 x_1})   \Big ] \\
&+ \frac{\mu}{2} \Big [ \big ( \nor \f^0_{,1}\nor^2 - 1  \big )^2 + \big (  \frac{1}{\sin ^2 x_1}\nor \f ^0_{,2}\nor^2 -1 \big)^2+ \big(  \nor \f^1_{,3}\nor^2 - 1 \big)^2 + \frac{2}{\sin^2 x_1} \big( \f^0_{,1}\cdot \f^0_{,2}\big )^2 + 2 \big( \f^0_{,1}\cdot \f^1_{,3}\big)^2 \\
&+  \frac{2}{\sin^2 x_1} \big( \f^0_{,2}\cdot \f^1_{,3}\big)^2 \Big ] + \frac{\lambda}{4} \Big(\nor \f^0_{,1}\nor^2 + \frac{1}{\sin ^2 x_1}\nor \f ^0_{,2}\nor^2 +\nor \f^1_{,3}\nor ^2  - 3  \Big )^2 \Big \} \sin x_1 dx.
\end{aligned}$$
Proceeding as in the cylindrical case, we consider the energy $\dis \w J$ defined by
$$\begin{aligned} 
&\w J(u) = \int_{\o}  \Big \{  \Big [ 2 \nor u_{,1}-a\nor ^2 + \frac{2}{\sin ^2 x_1} \nor u_{,2}-b\nor ^2  +   \nor\frac {1}{\sin^2 x_1} \; c+ u\nor^2   \Big ] \\
&+ \frac{\mu}{2} \Big [  \big(  \nor u\nor^2 - 1 \big)^2 + 2 \big( a\cdot u\big)^2 +  \frac{2}{\sin^2 x_1} \big( b\cdot u\big)^2 \Big ] + \frac{\lambda}{4} \Big(\nor a\nor^2 + \frac{1}{\sin ^2 x_1}\nor b\nor^2 +\nor u\nor ^2  - 3  \Big )^2 \Big \} \sin x_1 dx.
\end{aligned}$$
We notice that the $\dis \f^1_{,3}$ minimizing $J^0$ is the minimizer of $\w J$ with $\dis a= \f^0_{,1}$, $\dis b= \f^0_{,2}$ and $\dis c= \f^0_{,22}$. Next, considering the corresponding Euler-Lagrange equations, we obtain that the minimizer is a solution of the following partial differential equations
\begin{displaymath}
\begin{cases}
[\lambda (\nor a \nor ^2+ \frac{1}{\sin ^2 x_1} \nor b \nor ^2 +\nor u \nor^2 -3)+ 2 \mu(\nor u \nor ^2-1)+ 2 ]u+2\mu[(u,a)a+\frac{1}{\sin ^2 x_1} (u,b)b]\\
-4 \w\Delta u +\frac{2}{\sin x_1}c =0\; \text{in}\; \o\\
u=\w A a_3\; \text{on} \; \partial \o,
\end{cases}
\end{displaymath} 
where
\beq
 \w\Delta u= u_{,11} +\frac{1}{\sin^2 x_1}u_{,22},
\eeq
with $x_1$ being the latitudinal coordinate.

\paragraph{Acknowledgment.} I wish to thank Professor Herv\'e Le Dret for many useful discussions concerning the
subject of this paper.


\begin{thebibliography}{10}

\bibitem{ABP}
{\sc E. Acerbi, G. Buttazzo and D. Percivale}, A variational definition for the strain
energy of an elastic string, \emph{ J. Elasticity} 25, 137-148, 1991.

\bibitem{bZZ}
{\sc J. F. Babadjian, E. Zappale and H. Zorgati}, 
Dimensional reduction for energies with linear growth involving the bending moment, \emph{J. Math. Pures Appl.} (9) 90, No. 6, 520--549, 2008.

\bibitem{Ber}
{\sc M. Bernadou}, Finite Element Methods for Thin Shell Problems, \emph{John Wiley \& Sons Ltd.; Paris: Masson}, 1996.

\bibitem{bhajam}
{\sc K. Bhattacharya and R. D. James,}
A theory of thin films of martensitic materials with applications to microactuators,
\emph{J. Mech. Phys. Solids}, 47, 531--576, 1999.

\bibitem{BFF}
{\sc A. Braides, I. Fonseca and G. Francfort}, { 3D-2D asymptotic analysis for inhomogeneous thin films},  \emph{Indiana Univ. Math. J.} 49, No. 4, 1367-1404, 2000.

\bibitem{CR1}
{\sc G. Castineira and \'A. Rodriguez-Ar\'os}, Linear viscoelastic shells: an asymptotic approach, \emph{Asympt. Anal.} 107, No. 3--4, 169--201, 2018.


\bibitem{CR2}
{\sc G. Castineira and \'A. Rodriguez-Ar\'os}, On the justification of viscoelastic elliptic membrane shell equations, \emph{J. Elasticity},
130, No. 1, 85-113, 2018. 


\bibitem{chb}
{\sc D. Chapelle and K. J. Bathe,} The Finite Element Analysis of Shells. 2nd ed. Computational Fluid and Solid Mechanics, \emph{Berlin: Springer}, 410 p. 2011.

\bibitem{CDP}
{\sc X. Chen, H. Dai and E. Pruchnicki}, 
On a consistent rod theory for a linearized anisotropic elastic material. I: Asymptotic reduction method,
\emph{Math. Mech. Solids}, 26, No. 2, 217--229, 2021.


\bibitem{cia0}
{\sc P. G. Ciarlet,}
\emph{A justification of the Von K\'arm\'an equation.}
Arch. Rational Mech. Anal., 73: 349-189, 1980.

\bibitem{cia2}
P. G. \textsc{Ciarlet,}
\emph{Mathematical elasticity. Vol I: Three dimensional elasticity.}
North Holland, Amsterdam, 1988.

\bibitem{cia31}
P. G. \textsc{Ciarlet,}
\emph{Plates and Junctions in Elastic Multi-structures: An Asymptotic Analysis.}
Springer-Verlag, Berlin, 1990.


\bibitem{cia3}
P. G. \textsc{Ciarlet,}
\emph{Mathematical elasticity. Vol II: Theory of Plates.}
North Holland Publishing Co., Amsterdam, 1997.

\bibitem{cia4}
P. G. \textsc{Ciarlet,}
\emph{Mathematical elasticity. Vol III: Theory of Shells.}
North Holland Publishing Co., Amsterdam, 2000.

\bibitem{ciades}
P. G. \textsc{Ciarlet} and P. \textsc{Destuynder}
A justification of the two-dimensional linear plate model,
\emph{J. M\'ecanique,} 18(2): 315-344, 1979.


\bibitem{cialod}
{\sc P. G. Ciarlet and V. Lods},
Analyse asymptotique des coques lin\'eairement \'elastiques. I. Coques membranaires, \emph{C. R. Acad. Sci., Paris, S\'er. I} 319, No. 1, 863--868, 1994.

\bibitem{cialodmia}
{\sc P. G. Ciarlet, V. Lods and B. Miara},
Analyse asymptotique des coques lin\'eairement \'elastiques. I. Coques en flexion, \emph{C. R. Acad. Sci., Paris, S\'er. I} 319, No. 1, 95--100, 1994.

\bibitem{DAL}
{\sc G. Dal Maso}, { An Introduction to $\Gamma$-Convergence}, \emph{Birkhauser}, Boston, 1993.

\bibitem{DeG}
{\sc E. De Giorgi}, Sulla convergenza di alcune successioni di integrali del tipo dell?area.
\emph{Rend. Mat.}, (IV), 8  277-294, 1975.

\bibitem{DF}
{\sc E. De Giorgi and T. Franzoni}, Su un tipo di convergenza variazionale. \emph{Atti. Accad.
Naz. Lincei}, 58, 842-850, 1975.


\bibitem{FF}
{\sc I. Fonseca and G. Francfort}, { 3D-2D asymptotic analysis of an
	optimal design problem for thin films}, \emph{ journal fur die reine und angewandte mathematik}, vol. 1998, no. 505, pp. 173-202, 1998.

\bibitem{foxraosim}
D. D. \textsc{Fox}, A. \textsc{Raoult} and J.C. \textsc{Simo,}
A justification of nonlinear properly invariant plate theories,
\emph{Arch. Rational Mech. Anal.}, 124(2):157-199, 1993.


\bibitem{FD}
{\sc K. O. Friedrichs and R. F. Dressler},
A boundary-layer theory for elastic plates, \emph{Commun. Pure Appl. Math.} 14, 1-33, 1961.




\bibitem{Gol}
{\sc A. L. Goldenveizer,} Derivation of an approximate theory of shells by means of asymptotic integration of the equations of the theory of elasticity, \emph{J. Appl. Math. Mech.}, vol.27, no. 4, pp. 593-608, 1963.

\bibitem{IZ}
{\sc R. Ignat and H. Zorgati}, Dimension reduction and optimality of the uniform state in a phase-field-crystal model involving a higher-order functional, \emph{J. Nonlinear Sci.} 30, No. 1, 261--282, 2020.

\bibitem{ledmeu}
{\sc H. Le Dret and N. Meunier,}
Modeling heterogeneous martensitic wires, \emph{Math. Models Meth. Appl. Sci.}
Vol. 15 (03), 375-406, 2005.

\bibitem{ledmeu2}
{\sc H. Le Dret and N. Meunier,}
Heterogeneous wires made of martensitic materials,
\emph{C. R. Acad. Sci. Paris, Ser. I,} 337, 143--147, 2003.

\bibitem{ledrao}
{\sc H. Le Dret and A. Raoult,}
The nonlinear membrane model as variational limit of nonlinear three-dimensional elasticity,
	\emph{J. Math. Pures Appl.}, 75, 551--580, 1995.

\bibitem{ledrao2}
{\sc H. Le Dret and A. Raoult,}
The membrane shell model in nonlinear elasticity: A variational asymptotic derivation,
\emph{J. Nonlinear Sci.}, 6, 59--84, 1996.

\bibitem{ledzor}
{\sc H. Le Dret and H. Zorgati,}
Films courb\'es minces martensitiques,
\emph{C. R. Acad. Sci. Paris, Ser. I}, 339, 65--69, 2004.

\bibitem{ledzor2}
{\sc H. Le Dret and H. Zorgati,}
Asymptotic modeling of thin curved martensitic films,
\emph{ Asympt. Anal. ,} 48, 1-2, p. 141-171, 2006.

\bibitem{Meu}
{\sc N. Meunier}, 
Recursive derivation of one-dimensional models from three-dimensional nonlinear elasticity, \emph{Math. Mech. Solids}, 13, No. 2, 172-194 2008.



\bibitem{mia}
B. \textsc{Miara,}
\emph{Analyse asymptotique des coques membranaires non lin\'eairement \'elastiques.} C.R, Acad. Sci. Paris, 318, S\'erie I,p 689-694, 1994.


\bibitem{pan}
{\sc O. Pantz,}
\emph{Quelques probl\`emes de mod\'elisation en \'elasticit\'e non lin\'eaire,}
Doctoral dissertation, Universit\'e de Pierre et Marie Curie, Paris, 2001. 

\bibitem{rao}
A. \textsc{Raoult,}
\emph{Doctoral Dissertation.}
Universit\'e Pierre et Marie Curie, Paris, 1988.

\bibitem{SP1}
{\sc E. Sanchez-Palencia,} Statique et dynamique des coques minces. I. Cas de flexion pure non inhib\'ee, \emph{C. R. Acad. Sci. Paris}, 309, S\'erie I, pp. 411-417, 1989.

\bibitem{SP2}
{\sc E. Sanchez-Palencia,} Statique et dynamique des coques minces. II. Cas de flexion pure inhib\'ee, \emph{C. R. Acad. Sci. Paris}, 309, S\'erie I, pp. 531-537, 1989.

\bibitem{SP3}
{\sc E. Sanchez-Palencia,} Passage \`a la limite de l`\'elastiocit\'e tridimensionnelle \`a la th\'eorie asymptotique des coques minces, \emph{C. R. Acad. Sci. Paris}, 311, S\'erie II, pp. 909-916, 1990.

\bibitem{shu}
Y. C. {\sc Shu,}
Heterogeneous thin films of martensitic materials,
\emph{Arch. Rational Mech. Anal.}, 153, 39--90, (2000).

\bibitem{Tra}
{\sc K. Trabelsi},
Nonlinearly elastic thin plate models for a class of Ogden materials. I: The membrane model, \emph{Anal. Appl.}, Singap. 3, No. 2, 195--221 2005.

\bibitem{Zor1}
{\sc H. Zorgati},
A $\Gamma$-convergence result for thin curved films bonded to a fixed substrate with a noninterpenetration constraint,
\emph{Chin. Ann. Math.}, Ser. B 27, No. 6, 615-636 2006.

\bibitem{Zor2}
{\sc H. Zorgati},
Modeling thin curved ferromagnetic films,
\emph{Anal. Appl.}, Singap.	3, No. 4, 373-396 2005.

\bibitem{Zor3}
{\sc H. Zorgati},
Films courb\'es minces ferromagn\'etiques,
\emph{C. R., Math., Acad. Sci. Paris}, 340, No. 1, 81-86 2005.


\end{thebibliography}
\end{document}